\documentstyle{amsppt}
\pagewidth{5.11in} \pageheight{7.901in} \NoRunningHeads
\magnification = 1200

\topmatter
\title Monotonicity Formulae and Holomorphicity of Harmonic Maps between K\"ahler manifolds\endtitle
\author Yuxin Dong \endauthor
\thanks {*Supported by NSFC grant No 10971029, and NSFC-NSF grant No 1081112053}
\endthanks
\abstract {In this paper, we introduce the stress-energy tensors
of the partial energies $ E^{\prime}(f)$ and $E^{\prime
\prime}(f)$ of maps between K\"ahler manifolds. Assuming the
domain manifolds poss some special exhaustion functions, we use
these stress-energy tensors to establish some monotonicity
formulae of the partial energies of pluriharmonic maps into any
K\"ahler manifolds and harmonic maps into K\"ahler manifolds with
strongly semi-negative curvature respectively. These monotonicity
inequalities enable us to derive some holomorphicity and Liouville
type results for these pluriharmonic maps and harmonic maps. We
also use the stress-energy tensors to investigate the holomorphic
extension problem of $CR$ maps.}
\endabstract
\endsubjclass
\endkeywords
\endtopmatter
\document
\heading{\bf Introduction}
\endheading
\vskip 0.3 true cm

In 1980, Baird and Eells [BE] introduced the stress-energy tensor
for maps between Riemannian manifolds, which unifies various
results on harmonic maps. Following [BE], Sealey [Se] introduced
the stress-energy tensor for $p$-forms with values in vector
bundles and established some vanishing theorems for harmonic
$p$-forms. Since then, the stress-energy tensors have become a
useful tool for investigating the energy behavior of vector bundle
valued $p$-forms in various problems. Recently the authors in [DW]
presented a unified method to establish monotonicity formulae for
$p$-forms with values in vector bundles by means of the
stress-energy tensors of various energy functionals in geometry
and physics. Since the stress-energy tensors are $2$-tensor
fields, we may get $1$-forms by contracting them with vector
fields. The divergence of these $1$-forms then leads to a
fundamental integral formula, which is naturally linked to
conservation laws. The integral formula turns out to be a useful
tool for establishing monotonicity formulae of the energies,
provided that the $p$-forms satisfy the conservation laws and the
radial curvatures of the domain manifolds satisfy some pinching
conditions. Besides their possible applications in regularity
problems, these monotonicity formulae enable us not only to deduce
some vanishing theorems for $p$-forms under suitable growth
conditions, but also to investigate constant Dirichlet boundary
value problems for 1-forms. In [DW], the authors mainly used the
distance function of a complete Riemannian manifold to construct
the vector field in the integral formula. The Hessian of the
distance function appears naturally in the integral formula.
Consequently they used Hessian comparison theorems and coarea
formula to obtain their results. For the purposes of this paper,
we will consider more general exhaustion functions on domain
manifolds to construct vector fields in applying the integral
formula. Assuming that the domain manifolds poss some suitable
exhaustion functions, we may also establish some monotonicity
formulae for the $p$-forms which satisfy the conservation laws
(see Proposition 1.2 in \S 1).

In [Si1,2], Siu studied the holomorphicity of harmonic maps from
compact K\"ahler manifolds into compact K\"ahler manifolds with
strongly negative curvature or compact quotients of irreducible
symmetric bounded domains. The basic discovery of Siu was a
$\partial \overline{\partial }$-Bochner formula for harmonic maps
which does not involve the Ricci curvature tensor of the domains
(this is where K\"ahlerianity of the domains enters). Using his
modified Bochner formula and integration by parts, he proved a
vanishing theorem which implies that the harmonic maps in question
are actually pluriharmonic and some curvature terms of the
pull-back complexified tangent bundles vanish. The vanishing
curvature terms, under the assumption of sufficiently high rank,
forces the maps to be either holomorphic or anti-holomorphic.
Later, Sampson [Sa] extended Siu's technique to treat harmonic
maps of compact K\"ahler manifolds to more general targets. In
particular, he showed that harmonic maps from compact K\"ahler
manifolds to Riemannian manifolds with nonpositive Hermitian
curvature are pluriharmonic too. In [Li] and [PRS], the authors
generalized Sampson's pluriharmonicity result to noncompact
setting by assuming some growth conditions on energy of the
harmonic maps. On the other hand, some authors have investigated
the holomorphicity or pluriharmonicity of stable harmonic maps
from compact K\"ahler manifolds (cf. [SY], [Ud], [OU]). We refer
the reader to [To] for other related progress not mentioned here.

In this paper, we investigate harmonic maps between complete
K\"ahler manifolds. When the domain K\"ahler manifold is complete
(noncompact), the idea of the $\partial \overline{\partial
}$-Bochner technique in [Si1,2] together with the integration by
parts does not work any more and the holomorphicity problem of the
harmonic map is largely unknown. Notice that for a smooth map
$f:M\rightarrow N$ between two K\"ahler manifolds, one may
introduce two $1$-forms $\sigma $ and $\tau$ with values in
$f^{-1}TN$ whose vanishing characterizes the holomorphicity and
anti-holomorphicity of the map respectively. Actually $ \sigma
=\overline{\partial }f+\partial \overline{f}$ and $\tau =\partial
f+ \overline{\partial }\overline{f}$ if $\sigma $ and $\tau $ are
complexified (see \S2). The $L^2$ energies of $\sigma $ and $\tau
$ are just the partial energies $E^{\prime \prime}(f)$ and  $
E^{\prime}(f)$ respectively. Therefore we have the stress-energy
tensors $S_\sigma $ and $S_\tau$ corresponding to $\sigma$ and
$\tau$ respectively. It is natural to attempt to apply $S_\sigma $
and $S_\tau$ to investigate the energy behavior of the partial
energies and obtain vanishing theorems for $\sigma$ and $\tau$.
For this purpose, we assume that the domain manifolds poss some
suitable exhaustion functions. The advantages of using more
general exhaustion functions instead of the distance functions in
establishing monotonicity and vanishing results are that one may
not only relax the curvature conditions on the domain manifolds
but also has more choices for constructing suitable vector fields
in the integral formula. Fortunately some classes of complex
manifolds poss the required exhaustion functions. As the results
of this method, we obtain the pluriharmonicity of a harmonic map,
the monotonicity formulae for partial energies of a pluriharmonic
map or a harmonic map, the holomorphicity and constancy of a
pluriharmonic map or a harmonic map, and the holomorphic
extensions of $CR$ maps, etc.

Our method is based on the formulae (1.12), (1.13) for $\sigma$
and $\tau$, and computing $div(S_\sigma )$ and $div(S_\tau )$.
There are two $1$-forms $\gamma$ and $\rho$ arising naturally in
$div(S_\sigma )$ and $div(S_\tau)$. Then we derive the divergence
formulae of $\gamma $ and $\rho $, which are Weitzenb\"ock-type
formulae involving only the square norm of the $(1,1)$-part of the
second fundamental form $\nabla df$ and the curvature of the
target manifold. Assuming the domain K\"ahler manifold posses some
exhaustion function, these two divergence formulae enable us to
prove that a harmonic map into a K\"ahler manifold with strongly
semi-negative curvature is pluriharmonic if either
$|\overline{\partial } f|^2$ or $|\partial f|^2$ satisfies some
non-integrability condition (see Theorem 3.6). It follows that if
one of the partial energies has growth order at most $2$ (with
respect to the exhaustion function), then the harmonic map is
pluriharmonic (see Corollary 3.7) . In this way, we generalize
Siu's pluriharmonicity result to the non-compact setting. While
the authors in [PRS] considered more general targets in their
pluriharmonicity result by assuming the nonintegrability condition
on the energy density, we only assume the nonintegrability
condition on one of the partial energy densities to derive the
pluriharmonicity.

Next we investigate the monotonicity and holomorphicity of
harmonic maps between K\"ahler manifolds. First, we show that if
$f:M\rightarrow N$ is pluriharmonic, then $\sigma $ and $\tau $
satisfy the conservation laws, that is, $divS_\sigma =divS_\tau
=0$. Assuming $M$ posses a special exhaustion function (see
(4.1),(4.2) and (4.3)), it turns out that the conditions of
Proposition 1.2 are satisfied in this case. Hence we are able to
establish the monotonicity formulae for the partial energies of
the pluriharmonic map (see Theorem 4.3). It follows from the known
comparison theorems that the distance functions of some complete
K\"ahler manifolds become special exhaustion functions if their
radial curvatures have some suitable upper bounds (see Lemma 4.6).
This leads to the monotonicity formulae of pluriharmonic maps from
these complete K\"ahler manifolds. Remarkably no curvature
conditions are assumed on the targets for these results on
pluriharmonic maps. When a harmonic map between two K\"ahler
manifolds is not pluriharmonic, $\sigma $ and $\tau $ don't
satisfy the conservation laws in general. Due to this
non-conservativity, we can not apply Proposition 1.2 directly to
$\sigma$ and $\tau$. However, if the target K\"ahler manifold has
strongly semi-negative curvature, the fundamental integral
formulae related to the stress-energy tensors are still
applicable, because both $div(S_\sigma )$ and $div(S_\tau )$
contracted with suitable vector fields have some non-negativity
(see Lemma 4.9). Therefore one may establish the monotonicity
formulae of harmonic maps from certain K\"ahler manifolds to
K\"ahler manifolds with strongly semi-negative curvature too (see
Theorem 4.10). Besides the global monotonicity formulae, we also
obtain some local monotonicity formulae for partial energies of
pluriharmonic maps into K\"ahler manifolds or harmonic maps into
K\"ahler manifolds with strongly semi-negative curvature. Here, by
"local" we mean that the monotonicity formulae hold either in a
neighborhood of a point or outside a compact subset. Notice that
the authors in [DW] assumed some curvature pinching conditions to
establish monotonicity formulae for general $p$-forms. However,
the special properties (2.12) of $\sigma$ and $\tau$ enable us not
only to establish the monotonicity formulae on domain K\"ahler
manifolds whose curvatures only have some upper bounds, but also
to deduce the monotonicity formulae outside a compact subset. All
these monotonicity formulae imply immediately the holomorpicity of
the harmonic maps or pluriharmonic maps under suitable growth
conditions on the partial energies. In particular, Liouville type
theorems follow from suitable growth condition on the energy of
the maps. We should mention that a somewhat related approach has
been used by other authors, see e.g. [Ta 1,2]. However, he used a
integral formula technique to estimate the energy of harmonic maps
between K\"ahler manifolds (not the partial energies) and was
forced to get only Liouville type results. Our method of using
stress energy tensors seems to be easily operational and can also
be used to simplify the arguments in [Ta1,2].

The classical Bochner theorem [Bo] asserts that if $f$ is a smooth
CR function on the smooth connected boundary $\partial D$ of a
bounded domain $D$ in $C^m$, then $f$ can be extended from
$\partial D$ to $\overline{D}$ so that $f$ is holomorphic in $D$.
The Bochner type holomorphic extension problem for maps between
K\"ahler manifolds becomes a much harder problem, which has been
studied by several authors. In [Si2], [Sh], [NS], [Wo] and [CL],
the authors took a harmonic map approach by using Siu's $\partial
\overline{\partial }-$Bochner formula. The basic procedure for
this problem is as follows: First, one may find a harmonic
extension $f$ of the boundary map by the existence result of
Hamilton [Ha] and Schoen [Sc]. Next, try to derive the
holomorphcity of the harmonic extension.

Notice that the harmonic extension $f$ satisfies the tangential
Cauchy-Riemann equation on $\partial D$ if and only if the 1-form
$\sigma$ annihilates any tangent vector in the holomorphic
distribution $H$ of $\partial D$, that is, $\sigma|_H=0$. In their
generalization of Karcher-Wood theorem about constant Dirichlet
boundary problem for harmonic maps [KW], the authors in [DW],
using the stress-energy tensor and its related integral formula,
proved that if a $1$-form with values in a vector bundle satisfies
the conservation law over a starlike smooth domain $\Omega $ in
certain Riemannian manifold and annihilates any tangent vector of
$\partial \Omega $, then the $1-$form vanishes on $\Omega $. We
show that this method can also be applied to investigate the
vanishing of $\sigma$, although it only annihilates vectors in a
subbundle of $T(\partial D)$. Here (2.12) plays an important role
too. Consequently we are able to give an alternative proof of the
result in [CL](see Proposition 6.6) and obtain some other
holomorphic extension results not included in [Si2], [Sh] and [CL]
(see \S 6 for details).

This paper is organized as follows. In \S 1, we recall some basic
notions and formulae, and then describe briefly the approach of
[DW], but in a slightly generalized way.  In \S 2, we introduce
the stress-energy tensor $ S_\sigma $ and $S_\tau$ corresponding
to $\sigma $ and $\tau$ respectively. The relationship among the
stress-energy tensor $S_f$ introduced by Baird-Eells [BE] and the
stress-energy tensors $S_\sigma$, $S_\tau$ are discussed. In \S 3,
we give some criteria for the pluriharmonicity of a harmonic map
between K\"ahler manifolds. \S 4 and \S 5 are devoted to the
monotonicity formulae and holomorphicity of a pluriharmonic map or
a harmonic map between K\"ahler manifolds. Finally, in \S6, we
investigate the holomorphic extension problem of a CR boundary
map.

\heading{\bf 1. Monotonicity formulae and vanishing results of
$p$-forms}
\endheading
\vskip 0.3 true cm

Let $(M,g)$ be a Riemannian manifold and $\xi :E\rightarrow M$ a
Riemannian vector bundle over $M$ with a metric compatible
connection $\nabla ^E$. Let $A^p(\xi )$ denote the space of smooth
$p$-forms on $M$ with values in the vector bundle $\xi
:E\rightarrow M$, that is, $A^p(\xi )=\Gamma (\Lambda
^pT^{*}M\otimes E)$. The exterior covariant differentiation
$d^\nabla :A^p(\xi )\rightarrow A^{p+1}(\xi )$ relative to the
connection $\nabla ^E$ is defined by (cf. [EL])
$$
\aligned
d^\nabla \omega
(X_1,...,X_{p+1})&=\sum_{i=1}^{p+1}(-1)^{i+1}\nabla
_{X_i}^E(\omega (X_1,...,\widehat{X}_i,...,X_{p+1})) \\
&+\sum_{i<j}(-1)^{i+j}\omega
([X_i,X_j],X_1,...,\widehat{X}_i,...,\widehat{X} _j,...,X_{p+1})
\endaligned\tag{1.1}
$$
where the symbols covered by $\widehat{}$ are omitted. Since the
Levi-Civita connection on $TM$ is torsion-free, we also have
$$
(d^\nabla \omega
)(X_1,...,X_{p+1})=\sum_{i=1}^{p+1}(-1)^{i+1}(\nabla _{X_i}\omega
)(X_1,...,\widehat{X}_i,...,X_{p+1}).\tag{1.2}
$$
The induced inner product on $\Lambda ^pT_x^{*}M\otimes E_x$ is
defined as follows:
$$
<\alpha ,\beta >=\sum_{i_1<\cdots <i_p}<\alpha
(e_{i_1},...,e_{i_p}),\beta (e_{i_1},...,e_{i_p})>_{E_x} \tag{1.3}
$$
where $\{e_1,...,e_m\}$ is an orthonormal basis of $T_xM$,
$\forall x\in M$. Relative to the Riemannian structures of $E$ and
$TM$, the codifferential operator $\delta ^\nabla :A^p(\xi
)\rightarrow A^{p-1}(\xi )$ is characterized as the adjoint of
$d^\nabla $ via the formula:
$$
\int_M<d^\nabla \omega ,\theta >dv_g=\int_M<\omega ,\delta ^\nabla
\theta>dv_g
$$
where $\omega \in A^{p-1}(\xi )$, $\theta \in A^p(\xi )$, one of
which has compact support. Then
$$
(\delta ^\nabla \theta )(X_1,...,X_{p-1})=-\sum_i(\nabla
_{e_i}\theta )(e_i,X_1,...,X_{p-1}).  \tag{1.4}
$$
Let $f:M\rightarrow N$ be a smooth map between two Riemannian
manifolds. The pull-back bundle $f^{-1}TN$ is endowed with the
pull-back Riemannian structure. The differential $df$ may be
regarded as an element of $A^1(f^{-1}TN)$. The energy of $f$ is
defined by
$$
E(f)=\frac 12\int_M|df|^2dv_g.  \tag{1.5}
$$
In [BE], Baird-Eells introduced the stress-energy tensor $S_f$
associated with $E(f)$ as follows
$$
S_f=\frac{|df|^2}2g-df\odot df  \tag{1.6}
$$
where $df\odot df\in \Gamma (T^{*}M\otimes T^{*}M)$ is a symmetric
tensor defined by
$$
(df\odot df)(X,Y)=<df(X),df(Y)>.
$$
Then they proved that a harmonic map satisfies the conservation
law, that is, $divS_f=0$.

In general, we may introduce the following energy functional for
$\omega \in A^p(\xi )$
$$
F(\omega )=\int_M|\omega |^2dv_g . \tag{1.7}
$$
The stress-energy tensor associated with $F$ is a symmetric
$2$-tensor field given by (cf. [Se], [Xi], [Ba]):
$$
S_\omega (X,Y)=\frac{|\omega |^2}2g(X,Y)-(\omega \odot \omega
)(X,Y) \tag{1.8}
$$
for any $X,Y\in TM$, where $\omega \odot \omega $ denotes a
$2$-tensor field defined by
$$
(\omega \odot \omega )(X,Y)=<i_X\omega ,i_Y\omega > . \tag{1.9}
$$
Here $i_X\omega \in A^{p-1}(\xi )$ denotes the interior product by
$X\in TM$ , that is,
$$
(i_X\omega )(Y_1,...,Y_{p-1})=\omega (X,Y_1,...,Y_{p-1})
\tag{1.10}
$$
for any $Y_l\in TM$, $1\leq l\leq p-1$.

For a $2$-tensor field $T\in \Gamma (T^{*}M\otimes T^{*}M)$, its
divergence $ divT\in \Gamma (T^{*}M)$ is defined by
$$
(divT)(X)=\sum_i(\nabla _{e_i}T)(e_i,X),\ \ \forall X\in TM
\tag{1.11}
$$
where $\{e_i\}$ is an orthonormal basis of $TM$.

\proclaim{Lemma 1.1} (cf. [Se], [Xi], [Ba]) $(divS_\omega
)(X)=<\delta ^\nabla \omega ,i_X\omega
>+<i_Xd^\nabla \omega ,\omega >$
\endproclaim

\definition{Definition 1.1} $\omega \in A^p(\xi )$ is said to satisfy the
conservation law if $S_\omega $ is divergence free, that is,
$divS_\omega \equiv 0$.
\enddefinition

For a vector field $X$ on $M$, we denote by $\theta _X$ its dual
one form, that is,
$$
\theta _X(Y)=g(X,Y),\quad \forall \,Y\in TM.
$$
The covariant derivative of $\theta _X$ gives a $2$-tensor field
$\nabla \theta _X:$
$$
(\nabla \theta _X)(Y,Z)=(\nabla _Z\theta _X)(Y)=g(\nabla
_ZX,Y),\quad \forall Y,Z\in TM.
$$
If $X=\nabla \psi $ is the gradient of some smooth function $\psi
$ on $M$, then $\theta _X=d\psi $ and $\nabla \theta _X=Hess(\psi
)$.

For any vector field $X$ on $M$, a direct computation yields (cf.
Lemma 2.4 of [DW]):
$$
div(i_XS_\omega)=<S_\omega ,\nabla\theta_X>+(divS_\omega
)(X).\tag{1.12}
$$
Let $D$ be any bounded domain of $M$ with $C^1-$boundary. By
(1.12) and using the divergence theorem, we immediately have the
following integral formula (see also [Xi], [DW]):
$$
\int_{\partial D}S_\omega (X,\nu )dv_{\partial D}=\int_D[<S_\omega
,\nabla \theta _X>+(divS_\omega )(X)]dv_g  \tag{1.13}
$$
where $\nu $ is the unit outward normal vector field along
$\partial D$. In particular, if $\omega $ satisfies the
conservation law, then
$$
\int_{\partial D}S_\omega (X,\nu )=\int_D<S_\omega ,\nabla \theta
_X>.\tag{1.14}
$$
From now on, we often omit the volume elements in integral
formulae for simplicity when the integral domains are clear.

\definition{Definition 1.2} A function $\Phi :M\rightarrow R$ is called an
exhaustion function for a manifold $M$ if for every $t\in R$ the
sublevel set $\{x\in M:\Phi (x)<t\}$ is relatively compact in $M$.
The sublevel set $\{x\in M:\Phi (x)<t\}$ will be denoted by
$B_\Phi (t)$.
\enddefinition

Now we assume that $\Phi $ is a Lipschitz continuous exhaustion
function for a Riemannian manifold $M$ satisfying the following
conditions:

(1.15) $\Phi \geq 0$ and $\alpha =\sup_{x\in M}|\nabla \Phi |^2$
is finite;

(1.16) $\Psi =\Phi ^2$ is of class $C^\infty $ and $\Psi $ has
only discrete critical points.
\newline
It is a known fact that if a Riemannian manifold posses an
exhaustion function $\Phi$ with the property (1.15), then it is
complete. Actually let $r$ be the distance function relative a
point $o\in M$ and $x_k$ a Cauchy sequence of $M$. The triangle
inequality implies that $r(x_k)\le C$ for some constant $C$. It
follows that $|\Phi(x_k)|\le |\Phi(o)|+\sqrt{\alpha }r(x_k)\le
|\Phi(o)|+\sqrt{\alpha }C$. By the properness of $\Phi$, we see
that the sequence $x_k\rightarrow x_0\in M$, that is, $M$ is
complete.

\proclaim{Proposition 1.2} Let $M$ be a Riemannian manifold and
let $\Phi $ be an exhaustion function with properties (1.15) and
(1.16). Suppose $\xi :E\rightarrow M$ is a Riemannian vector
bundle on $M$ and $\omega \in A^p(\xi )$ satisfies the
conservation law. If there exists a positive constant $\beta $
such that
$$
<S_\omega ,Hess(\Psi )>\geq \beta |\omega |^2  \tag{1.17}
$$
then
$$
\frac 1{\rho _1^\Lambda }\int_{B_\Phi (\rho _1)}|\omega |^2\leq
\frac 1{\rho _2^\Lambda }\int_{B_\Phi (\rho _2)}|\omega
|^2\tag{1.18}
$$
for any $0<\rho _1\leq \rho _2$, where $\Lambda =\beta /\alpha $.
Furthermore, if $ \int_{B_\Phi (t)}|\omega |^2=o(t^\Lambda)$ $($
as $t\rightarrow \infty )$, then $\omega =0$.
\endproclaim

\demo{Proof} Set $X=\frac 12\nabla (\Psi )=\Phi \nabla \Phi $.
Obviously $(\nabla \Psi )|_{\partial B_\Phi (t)}$ is an outward
normal vector field along $\partial B_\Phi (t)$  if $t>0$ is a
regular value of $\Psi $. By (1.8), we have
$$
\aligned
S_\omega (X,\nu )&=t<\nabla \Phi ,\nu >[\frac{|\omega
|^2}2-|i_\nu \omega |^2]
\\
&\leq  \frac{t  \sqrt{\alpha }|\omega |^2}2
\endaligned \tag{1.19}
$$
on $\partial B_\Phi (t)$. It follows from (1.14), (1.17) and
(1.19) that
$$
t\sqrt{\alpha }\int_{\partial B_\Phi (t)}\frac{|\omega |^2}2 \geq
\frac \beta 2\int_{B_\Phi (t)}|\omega |^2.\tag{1.20}
$$
By the co-area formula, we have
$$
\aligned \frac d{dt}\int_{B_\Phi (t)}|\omega |^2&=\frac
d{dt}\{\int_0^t\int_{\partial
B_\Phi (s)}\frac{|\omega |^2}{|\nabla \Phi |})ds \\
&=\int_{\partial B_\Phi (t)}\frac{|\omega |^2}{|\nabla \Phi |} \\
&\geq \frac 1{\sqrt{\alpha }}\int_{\partial B_\Phi (t)}|\omega
|^2.
\endaligned \tag{1.21}
$$
Dividing both sides of (1.20) by $\alpha $, we obtain from(1.20)
and (1.21) that
$$
t\frac d{dt}\int_{B_\Phi (t)}|\omega |^2\geq\Lambda \int_{B_\Phi
(t)}|\omega |^2
$$
where $\Lambda =\beta /\alpha $. Hence we get
$$
\frac d{dt}(\frac{\int_{B_\Phi (t)}|\omega |^2}{t^\Lambda })\geq
0. \tag{1.22}
$$
This proposition follows immediately from integrating (1.22) on
$[\rho _1,\rho _2]$. \qed
\enddemo

\remark{Remark 1.1}
\newline (a) A more general monotonicity inequality holds if one further assumes
either $|\nabla \Phi |=const.$ or $\frac{|\omega |^2}2-|i_\nu
\omega |^2\geq 0$ a.e. on $M$. In both cases, the first line of
(1.19) yields
$$
S_\omega (X,\nu )\leq t\sqrt{\alpha }(\frac{|\omega |^2}2-|i_\nu
\omega |^2)
$$
and thus
$$t\sqrt{\alpha }\int_{\partial B_\Phi (t)}|\omega
|^2-\beta \int_{B_\Phi (t)}|\omega |^2\geq 2t\sqrt{\alpha
}\int_{\partial B_\Phi (t)}|i_\nu \omega |^2.$$ Similar arguments
imply that
$$
\frac d{dt}(\frac{\int_{B_\Phi (t)}|\omega |^2}{t^\Lambda })\geq \frac{%
2t^{-\Lambda }}{\sqrt{\alpha }}\int_{\partial B_\Phi (t)}|i_\nu
\omega |^2
$$
Hence
$$
\frac 1{\rho _2^\Lambda }\int_{B_\Phi (\rho _2)}|\omega |^2-\frac
1{\rho
_1^\Lambda }\int_{B_\Phi (\rho _1)}|\omega |^2\geq \frac 2{\sqrt{\alpha }%
}\int_{\rho _1}^{\rho _2}\frac{\int_{\partial B_\Phi (t)}|i_\nu \omega |^2}{%
t^\Lambda }dt
$$
for any $0<\rho _1\leq \rho _2$.
\newline (b) The proof of
Proposition 1.2 shows that the monotonicity formula still holds if
one assumes $\int_{B_\Phi (t)}divS_\omega (X)\geq 0$ for $t>0$
instead of assuming $divS_\omega =0$. We will consider these
important non-conservative cases too.
\endremark

We will apply Proposition 1.2 to investigate the monotonicity and
holomorphicity of harmonic maps and pluriharmonic maps between
K\"ahler manifolds.

\heading{\bf 2. The stress-energy tensors of $\overline{\partial
}f$ and $\partial f $}
\endheading
\vskip 0.3 true cm

A Hermitian metric on a complex manifold $M$ is a Riemannian
metric $g$ such that $g(JX,JY)=g(X,Y)$, $\forall X,Y\in TM$, where
$J$ denotes the complex structure of $M$. We say that $(M,g)$ is
K\"ahler if $J$ is parallel with respect to the Levi-Civita
connection of $g$, that is, $\nabla J=0$.

We denote by $<\cdot ,\cdot >$ the (real) inner product of tensor
bundles of $M$ induced by $g$. The complex extension of the inner
product is still denoted by $<\cdot ,\cdot >$. Define the
Hermitian inner product $\ll \cdot ,\cdot \gg $ by
$$
\ll u,v\gg =<u,\overline{v}>.  \tag{2.1}
$$

Henceforth $(M^m,g)$ and $(N^n,h)$ will denote K\"ahler manifolds
of complex dimensions $m$ and $n$ respectively. Let
$f:M\rightarrow N$ be a smooth map from $M$ to $N$. The complex
structure of $M$ (resp. $N$) gives a decompositions of $TM^C$
(resp. $TN^C$) into tangent vectors of type $(1,0)$ and type
$(0,1)$. Then we have
$$
TM^C=T^{1,0}M\oplus T^{0,1}M,\text{\ }TN^C=T^{1,0}N\oplus
T^{0,1}N.\tag{2.2}
$$
By restricting and projecting the complexified differential
$df:TM^C\rightarrow TN^C$ to the subbundles in (2.2), we have the
following bundle maps (cf. also [Si2]):
$$
\aligned &\partial f:T^{1,0}M\rightarrow T^{1,0}N,\quad
\overline{\partial}f:T^{0,1}M\rightarrow T^{1,0}N \\
&\partial \overline{f}:T^{1,0}M\rightarrow T^{0,1}N, \quad
\overline{\partial f }:T^{0,1}M\rightarrow T^{0,1}N.
\endaligned\tag{2.3}
$$
The energy functional of maps is defined by $E(f)=\frac
12\int_M|df|^2$, where the energy density is
$$
\frac 12|df|^2=\frac
12\sum_{j=1}^m\{<df(e_j),df(e_j)>+<df(Je_j),df(Je_j)>\} \tag{2.4}
$$
in terms of an orthonormal basis $\{e_i,Je_i\}_{i=1,...,m}\in TM$.
Set
$$
\eta _j=\frac 1{\sqrt{2}}(e_j-iJe_j),\quad \eta
_{\overline{j}}=\overline{ \eta }_j=\frac
1{\sqrt{2}}(e_j+iJe_j).\tag{2.5}
$$
Then $\ll \eta _j,\eta _k\gg =\ll \eta _{\overline{j}},\eta
_{\overline{k} }\gg =\delta _{jk}$, that is $\{\eta _j\}_{j=1}^m$
(resp. $\{\eta _{ \overline{j}}\}_{j=1}^m$) is a unitary basis of
$T_x^{1,0}M$ (resp. $ T_x^{0,1}M$). The partial energy densities
of $f$ are defined by
$$
|\overline{\partial }f|^2=\sum_{j=1}^m\ll \overline{\partial }f(\overline{%
\eta }_j),\overline{\partial }f(\overline{\eta }_j)\gg ,\quad
|\partial f|^2=\sum_{j=1}^m\ll \partial f(\eta _j),\partial f(\eta
_j)\gg.\tag{2.6}
$$
A direct computation gives
$$
\aligned |\overline{\partial }f|^2=\frac
14\sum_{j=1}^m\{&<df(e_j),df(e_j)>+<df(Je_j),df(Je_j)> \\
&-2<df(Je_j),J^{\prime }df(e_j)>\}
\endaligned
$$
and
$$
\aligned |\partial f|^2=\frac
14\sum_{j=1}^m\{&<df(e_j),df(e_j)>+<df(Je_j),df(Je_j)>
\\
&+2<df(Je_j),J^{\prime }df(e_j)>\}.
\endaligned
$$
Hence we have
$$
E(f)=E^{\prime }(f)+E^{\prime \prime }(f)
$$
where $E^{\prime }(f)=\int_M|\partial f|^2$ and $E^{\prime \prime
}(f)=\int_M| \overline{\partial }f|^2$ are the partial energies of
$f$ respectively. The map $f:M\rightarrow N$ is called holomorphic
(resp. anti-holomorphic) if $$df\circ J=J\circ df\quad(\text{resp.
}df\circ J=-J^{\prime }\circ df)$$ which is equivalent to
$\overline{\partial }f=0$ (resp. $\partial f=0$).

For a smooth map $f:M\rightarrow N$, we introduce two $1-$forms
$\sigma $, $\tau \in A^1(f^{-1}TN)$ as follows:
$$
\sigma (X)=\frac{df(X)+J^{\prime }df(JX)}2  \tag{2.7}
$$
and
$$
\tau (X)=\frac{df(X)-J^{\prime }df(JX)}2  \tag{2.8}
$$
for any $X\in TM$. Then (2.7) and (2.8) yield
$$
\sigma (JX)=\frac{df(JX)-J^{\prime }df(X)}2=-J^{\prime }\sigma (X)
\tag{2.9}
$$
and
$$
\tau (JX)=\frac{df(JX)+J^{\prime }df(X)}2=J^{\prime }\tau (X).
\tag{2.10}
$$
By complex extension, we may regard $\sigma $ and $\tau $ as
sections of $T^{*}M^C\otimes f^{-1}TN^C$. By (2.9) and (2.10), we
have
$$
\aligned \sigma&:T^{0,1}M\rightarrow T^{1,0}N,\quad \sigma :
T^{1,0}M\rightarrow T^{0,1}N \\
\tau&:T^{1,0}M\rightarrow T^{1,0}N,\quad \tau :T^{0,1}M\rightarrow
T^{0,1}N \endaligned\tag{2.11}
$$
and
$$
\aligned
&<\sigma (JX),\sigma (JY)>=<\sigma (X),\sigma (Y)> \\
&<\tau (JX),\tau (JY)>=<\tau (X),\tau (Y)>.
\endaligned\tag{2.12}
$$
It turns out that (2.12) is important for studying the paritial
energies of a map. By (2.7) and (2.8), we also derive
$$
\aligned |\sigma |^2=&\sum_{j=1}^m[<\sigma (e_j),\sigma
(e_j)>+<\sigma (Je_j),\sigma
(Je_j)>] \\
=&\frac
12\sum_{j=1}^m[<df(e_j),df(e_j)>+<df(Je_j),df(Je_j)>\\
&-2<df(Je_j),J^{\prime }df(e_j)>]
\endaligned
$$
and
$$
\aligned
|\tau |^2=&\sum_{j=1}^m[<\tau (e_j),\tau (e_j)>+<\tau (Je_j),\tau (Je_j)>] \\
=&\frac12\sum_{j=1}^m[<df(e_j),df(e_j)>+<df(Je_j),df(Je_j)>\\
&+2<df(Je_j),J^{\prime }df(e_j)].
\endaligned
$$
Therefore
$$
\aligned &|\sigma |^2+|\tau |^2=|df|^2\\
&|\sigma |^2=2|\overline{\partial }f|^2,\quad |\tau |^2=2|\partial
f|^2. \endaligned\tag{2.13}
$$

Recall that the map $f$ is said to be {\sl harmonic} if it
satisfies the Euler-Lagrange equation of the energy functional
$E(f)$, that is,
$$
\sum_{A=1}^{2m}(\nabla df)(e_A,e_A)=2\sum_{j=1}^m(\nabla df)(\eta
_j,\eta _{ \overline{j}})=0\tag{2.14}
$$
where $\nabla df$ denotes the second fundamental form of $f$ ([cf.
[EL]).

\proclaim{Lemma 2.1}
If $f:M\rightarrow N$ is a harmonic map between K\"ahler manifolds, then $%
\delta ^\nabla \sigma =\delta ^\nabla \tau =0$.
\endproclaim

\demo{Proof} Let $\{e_A\}_{A=1,...,2m}=\{e_i,Je_i\}_{i=1,...,m}$
be an orthonormal frame field around $p\in M$ such that $(\nabla
_{e_A}e_B)_p=0$. Since $f$ is harmonic, we have
$$
\aligned \delta ^\nabla \sigma &=-\sum_{A=1}^{2m}(\nabla
_{e_A}\sigma)(e_A)\\
&=-\frac 12\sum_{A=1}^{2m}\nabla _{e_A}[df(e_A)+J^{\prime }df(Je_A)] \\
&=-\frac 12J^{\prime }(\nabla df)(Je_A,e_A) \\
&=-\frac{J^{\prime }}2\sum_{i=1}^m[(\nabla df)(Je_i,e_i)-(\nabla
df)(e_i,Je_i)] \\
&=0. \endaligned
$$
Likewise we have $\delta ^\nabla \tau \equiv 0$. \qed
\enddemo

By definition, the stress-energy tensors of $\sigma $ and $\tau $
are given as follows:
$$
\aligned
S_\sigma (X,Y)&=\frac{|\sigma |^2}2g(X,Y)-\sigma \odot \sigma (X,Y) \\
S_\tau (X,Y)&=\frac{|\tau |^2}2g(X,Y)-\tau \odot \tau (X,Y).
\endaligned \tag{2.15}
$$
Using (2.12), we obtain $S_\sigma (JX,JY)=S_\sigma (X,Y)$ and $
S_\tau (JX,JY)=S_\tau (X,Y)$. This shows that $S_\sigma $ and
$S_\tau $ are $(1,1)-$type tensor fields, according to the
decomposition $T^{*}M^C\otimes T^{*}M^C=T^{*(2,0)}M\oplus
T^{*(1,1)}M\oplus T^{*(0,2)}M$.

We hope to find the relationship among the three stress-energy
tensors $S_f$ , $S_\sigma $ and $S_\tau $. From (1.6), (2.13) and
(2.15), we get
$$
S_f-[S_\sigma +S_\tau ]=\frac 12[(dfJ)\odot (dfJ)-df\odot
df].\tag{2.16}
$$
Write $\Psi _f=\frac 12[(df\circ J)\odot (df\circ J)-df\odot df]$.
For any two vector $Z,W\in T^{1,0}M$, we obtain
$$
\aligned \Psi _f(Z,\overline{W})&=\frac 12[h(idf(Z),-idf(\overline{W}))-h(df(Z),df(\overline{W}))] \\
&=0.\endaligned
$$
This implies
$$
S_f^{(1,1)}=S_\sigma +S_\tau, \quad \Psi
_f=S_f^{(2,0)}+S_f^{(0,2)}. \tag{2.17}
$$

\definition{Definition 2.1} A map $f:M\rightarrow N$ is called
{\sl pluriconformal} if $J$ is isometry w.r.t. $df\odot df$, that
is, $\Psi _f\equiv 0$.
\enddefinition

For any $Z=X-iJX$, $W=Y-iJY\in T^{1,0}M$, we have
$$
\aligned (f^{*}h)(Z,W)&=(df\odot df)(X,Y)-(df\odot df)(JX,JY) \\
& -i[(df\odot df)(X,JY)+(df\odot df)(JX,Y)].\endaligned
$$
Therefore $f$ is pluriconformal if and only if
$(f^{*}h)^{(2,0)}=0$ or equivalently $(f^{*}h)^{(0,2)}=0$. Clearly
$\pm $holomorphic maps are pluriconformal. However the converse is
not true in general. By (2.16), $ S_f=S_\sigma +S_\tau $ if and
only if $f$ is pluriconformal. Note also that the notion of
pluriconformal maps may be defined for maps from a K\"ahler
manifold to a Riemannian manifold, that is, the target manifold is
not necessarily a K\"ahler manifold.

Let $f:M\rightarrow N$ be a smooth map from a K\"ahler manifold.
The complexified second fundamental form $\nabla df$ of $f$ , as a
section of $ T^{*}M^C\otimes T^{*}M^C\otimes f^{-1}TN^C$, may be
decomposed as follows:
$$
\nabla df=(\nabla df)^{(2,0)}+(\nabla df)^{(1,1)}+(\nabla
df)^{(0,2)}. \tag{2.18}
$$

\definition{Definition 2.2}
A smooth map $f:M\rightarrow N$ from a K\"ahler manifold is called
{\sl pluriharmonic} if $(\nabla df)^{(1,1)}\equiv 0$, that is,
$$
(\nabla df)(Z,\overline{W})=0  \tag{2.19}
$$
for any $Z,W\in T^{(1,0)}M$.
\enddefinition

We should mention that the notion of pluriharmonic maps is
well-defined for any smooth map $\varphi :M\rightarrow
\widetilde{N}$ from a K\"ahler manifold to a Riemannian manifold
$\widetilde{N}$. When $\widetilde{N}=R$, $(\nabla d\varphi
)^{(1,1)}$ will be called the complex Hessian of $\varphi $ and
denoted by $H(\varphi )$. Note that any pluriharmonic map is
automatically harmonic, and any $\pm $holomorphic map between
K\"ahler manifolds is pluriharmonic too. Clearly a smooth map
between two K\"ahler manifolds is pluriharmonic if and only if its
restriction to every holomorphic curve in $M$ is harmonic (cf.
[Ra] for more general cases). The notion of pluriharmonic maps
lies between those of harmonic and $\pm $holomorphic maps. There
is no difference between harmonic and pluriharmonic for the case
$\dim _CM=1$.

For any $Z=X_1-iJX_1$ , $W=X_2-iJX_2$ $\in T^{(1,0)}M$, we have
$$
\aligned (\nabla df)(Z,\overline{W})=&(\nabla df)(X_1,X_2)+(\nabla
df)(JX_1,JX_2) \\
&+i[(\nabla df)(X_1,JX_2)-(\nabla df)(JX_1,X_2)].\endaligned
$$
This shows that (2.19) is equivalent to
$$
(\nabla df)(X_1,X_2)+(\nabla df)(JX_1,JX_2)=0. \tag{2.20}
$$
and
$$
(\nabla df)(X_1,JX_2)-(\nabla df)(JX_1,X_2)=0. \tag{2.21}
$$
for any $X_1,X_2\in TM$. Define a $2-$form $\alpha \in
A^2(f^{-1}TN)$ by
$$
\alpha (X_1,X_2)=(\nabla df)(X_1,JX_2)-(\nabla df)(JX_1,X_2).
\tag{2.22}
$$
for any $\,X_1,X_2\in TM$. By (2.22), we obtain $\alpha
(JX_1,JX_2)=\alpha (X_1,X_2)$. Writing $X_2=-J \widetilde{X}_2$ in
(2.21), we see that (2.20) and (2.21) are actually equivalent.
Hence we have proved the following:

\proclaim{Lemma 2.2} A map $f:M\rightarrow N$ is pluriharmonic$\
$if and only if $f$ satisfies (2.20), or equivalently $f$
satisfies (2.21), that is, $\alpha \equiv 0$.
\endproclaim

Let $f:M\rightarrow N$ be a smooth map between K\"ahler manifolds
with $y=f(x)$, $x\in M$. We may choose normal orthnormal frame
fields $\{e_j,Je_j\}_{j=1}^m$ and $\{\widetilde{e}_\alpha
,J^{\prime }\widetilde{e}_\alpha \}_{\alpha =1}^n$ around $x$ and
$y$ respectively. Set
$$
\aligned \eta _j&=\frac 1{\sqrt{2}}(e_j-iJe_j),\quad \eta
_{\overline{j}}=\overline{
\eta }_j=\frac 1{\sqrt{2}}(e_j+iJe_j) \\
\widetilde{\eta }_\alpha &=\frac 1{\sqrt{2}}(\widetilde{e}_\alpha
-iJ^{\prime }\widetilde{e}_\alpha ),\quad \eta _{\overline{\alpha
}}=\overline{\eta } _\alpha =\frac
1{\sqrt{2}}(\widetilde{e}_\alpha +iJ\widetilde{e}_\alpha).
\endaligned \tag{2.23}
$$
Put $f_j=df(\eta _j)$ and $f_{\overline{j}}=df(\eta
_{\overline{j}})$. From (2.2), we have
$$
\aligned
f_j&=f_j^{(1,0)}+f_j^{(0,1)}=(\partial f)(\eta _j)+(\partial \overline{f}%
)(\eta _j) \\
f_{\overline{j}}&=f_{\overline{j}}^{(1,0)}+f_{\overline{j}}^{(0,1)}=(%
\overline{\partial }f)(\eta _{\overline{j}})+(\overline{\partial }\overline{f%
})(\eta _{\overline{j}}).
\endaligned  \tag{2.24}
$$
Using the frame field $\{\widetilde{\eta }_\alpha ,\widetilde{\eta
}_{ \overline{\alpha }}\}_{\alpha=1,...,n}$, we may write $f_j$
and $f_{\overline{j}}$ as follows:
$$
\aligned f_j&=\sum_\alpha (f_j^\alpha \widetilde{\eta }_\alpha
+f_j^{\overline{\alpha }
}\widetilde{\eta }_{\overline{\alpha }}) \\
f_{\overline{j}}&=\sum_\alpha (f_{\overline{j}}^\alpha
\widetilde{\eta } _\alpha +f_{\overline{j}}^{\overline{\alpha
}}\widetilde{\eta }_{\overline{ \alpha }}).
\endaligned \tag{2.25}
$$
Let $\{\theta ^j,\theta ^{\overline{j}}\}_{j=1}^m$ and
$\{\widetilde{\theta } ^\alpha ,\widetilde{\theta
}^{\overline{\alpha }}\}_{\alpha =1}^n$ be the dual frame fields
of $\{\eta _j,\eta _{\overline{j}}\}$ and $\{\widetilde{\eta
}_\alpha ,\widetilde{\eta }_{\overline{\alpha }}\}$ respectively.
Hence we may express the complexified second fundamental form
$\nabla df$ as follows
$$
\aligned \nabla df=&\sum_{k,j,\alpha }\{f_{kj}^\alpha \theta
^k\otimes \theta ^j\otimes \widetilde{\eta }_\alpha
+f_{kj}^{\overline{\alpha }}\theta
^k\otimes \theta ^j\otimes \widetilde{\eta }_{\overline{\alpha }} \\
&+f_{\overline{k}j}^\alpha \theta ^{\overline{k}}\otimes \theta
^j\otimes \widetilde{\eta }_\alpha
+f_{\overline{k}j}^{\overline{\alpha }}\theta ^{
\overline{k}}\otimes \theta ^j\otimes \widetilde{\eta
}_{\overline{\alpha }}
\\
&+f_{k\overline{j}}^\alpha \theta ^k\otimes \theta
^{\overline{j}}\otimes \widetilde{\eta }_\alpha
+f_{k\overline{j}}^{\overline{\alpha }}\theta ^k\otimes \theta
^{\overline{j}}\otimes \widetilde{\eta }_{\overline{\alpha }
} \\
&+f_{\overline{k}\overline{j}}^\alpha \theta
^{\overline{k}}\otimes \theta ^{ \overline{j}}\otimes
\widetilde{\eta }_\alpha +f_{\overline{k}\overline{j}}^{
\overline{\alpha }}\theta ^{\overline{k}}\otimes \theta
^{\overline{j} }\otimes \widetilde{\eta }_{\overline{\alpha }}\}.
\endaligned\tag{2.26}
$$
Then
$$
\aligned f_{kj}^\alpha &=f_{jk}^\alpha ,\quad
f_{k\overline{j}}^\alpha
=f_{\overline{j}k}^\alpha \\
\overline{f_{kj}^\alpha
}&=f_{\overline{k}\overline{j}}^{\overline{\alpha } },\quad
\overline{f_{k\overline{j}}^\alpha
}=f_{\overline{k}j}^{\overline{\alpha }}.
\endaligned\tag{2.27}
$$
and
$$
\aligned (\nabla df)^{(1,1)}&=\sum_{k,j,\alpha
}\{f_{k\overline{j}}^\alpha \theta ^k\otimes \theta
^{\overline{j}}\otimes \widetilde{\eta }_\alpha +f_{k
\overline{j}}^{\overline{\alpha }}\theta ^k\otimes \theta
^{\overline{j}
}\otimes \widetilde{\eta }_{\overline{\alpha }} \\
&+f_{\overline{k}j}^\alpha \theta ^{\overline{k}}\otimes \theta
^j\otimes \widetilde{\eta }_\alpha
+f_{\overline{k}j}^{\overline{\alpha }}\theta ^{
\overline{k}}\otimes \theta ^j\otimes \widetilde{\eta
}_{\overline{\alpha } }\}.
\endaligned\tag{2.28}
$$
From (2.28), we get
$$
|(\nabla df)^{(1,1)}|^2=4\sum_{k,j,\alpha
}|f_{k\overline{j}}^\alpha |^2.\tag{2.29}
$$
It follows that $f$ is pluriharmonic if and only if
$f_{k\overline{j} }^\alpha =0$, $1\leq k,j\leq m$, $1\leq \alpha
\leq n$.

\heading{\bf 3. Pluriharmonicity of harmonic maps}
\endheading
\vskip 0.3 true cm First we recall some curvature conditions
introduced by Y.T. Siu [Si1] (cf. also [Sa], [To], [OU]). Let
$(M^m,g)$ be a K\"ahler manifold of complex dimension $m$. The
curvature tensor $R$ of $M$ is defined by
$$
R(X,Y)Z=\nabla _X\nabla _YZ-\nabla _Y\nabla _XZ-\nabla
_{[X,Y]}Z\,,\ \ \forall \,X,Y,Z\in TM.  \tag{3.1}
$$
We denote by $Q$ the curvature operator defined by $R$
$$
<Q(X\wedge Y),Z\wedge W>=<R(X,Y)W,Z>,\ \ \forall \,X,Y,Z,W\in TM.
\tag{3.2}
$$
The complex extension of $Q$ to $\wedge ^2TM^C$ is also denoted by
$Q$. By (2.1), we have
$$
\ll Q(X\wedge Y),Z\wedge W\gg =<Q(X\wedge Y),\overline{Z\wedge
W}>,\ \ \forall \,X,Y,Z,W\in TM^C.  \tag{3.3}
$$
The K\"ahler identity of $M$ yields
$$
Q|\wedge ^{(2,0)}TM=Q|\wedge ^{(0,2)}TM=0. \tag{3.4}
$$
Set
$$
Q^{(1,1)}=Q:\wedge ^{(1,1)}TM\rightarrow \wedge ^{(1,1)}TM.
\tag{3.5}
$$
\definition{Definition 3.1 ([Si1])} The curvature tensor of $(M,g)$ is said to
be strongly negative (resp. strongly semi-negative) if
$$
\ll Q^{(1,1)}(\xi ),\xi \gg =<Q^{(1,1)}(\xi ),\overline{\xi
}>\!\text{ }<0\ \text{ (resp.\negthinspace }\leq 0\text{)}
\tag{3.6}
$$
for any $\xi =(Z\wedge W)^{(1,1)}\neq 0$, $\,Z,W\in TM^C$.
\enddefinition

Writing $Z=Z^{(1,0)}+Z^{(0,1)},\ W=W^{(1,0)}+W^{(0,1)}$ in
Definition 3.1, we get
$$
\xi =Z^{(1,0)}\wedge W^{(0,1)}-W^{(1,0)}\wedge Z^{(0,1)}.
\tag{3.7}
$$
and thus
$$
\aligned \ll Q(Z\wedge W),Z\wedge W\gg =\ll &Q(Z^{(1,0)}\wedge W^{(0,1)}-W^{(1,0)}\wedge Z^{(0,1)}, \\
&Z^{(1,0)}\wedge W^{(0,1)}-W^{(1,0)}\wedge Z^{(0,1)}\gg.
\endaligned\tag{3.8}
$$
In [Si1,Si2], Y.T. Siu showed the following result by his
$\partial \overline{\partial }-$Bochner formula:

\proclaim{Proposition 3.1 ([Si1])} Let $f:M\rightarrow N$ be a
harmonic map from a compact K\"ahler manifold into a K\"ahler
manifold with strongly semi-negative curvature. Then $f$ is a
pluriharmonic map and
$$
\ll \widetilde{Q}(f_j\wedge f_k),f_j\wedge f_k\gg =0,\quad
j,k=1,...,m \tag{3.9}
$$
where $\widetilde{Q}$ denotes the curvature operator of $N$.
\endproclaim

When $N$ is a K\"ahler manifold with strongly negative curvature
or an irreducible symmetric bounded domain, Siu derived the
holomorphicity of $f$ under further rank condition on $df$ (cf.
also Lemma 5.9 below).

We will extend Proposition 3.1 to the complete noncompact case.
Suppose $f:M^m\rightarrow N^n$ is a harmonic map between K\"ahler
manifold. Let $\{\eta _j,\eta _{\overline{j}}\}_{j=1}^m$ and
$\{\widetilde{\eta }_\alpha , \widetilde{\eta }_{\overline{\alpha
}}\}_{\alpha =1}^n$ be the normal unitary frame fields defined in
(2.27) and let $\{\theta ^j,\theta ^{ \overline{j}}\}_{j=1}^m$ and
$\{\widetilde{\theta }^\alpha ,\widetilde{ \theta
}^{\overline{\alpha }}\}_{\alpha =1}^n$ be their dual frame fields
respectively. By the definition of $\sigma $, we obtain
$$
\sigma (\eta _j)=[df(\eta_j)]^{(0,1)}=\sum_\alpha
f_j^{\overline{\alpha }}\widetilde{\eta }_{\overline{\alpha }}.
\tag{3.10}
$$
The complex conjugate of $\sigma(\eta_j)$ gives
$$
\sigma (\eta _{\overline{j}})=\sum_\alpha f_{\overline{j}}^\alpha
\widetilde{\eta }_\alpha. \tag{3.11}
$$
Hence
$$
\aligned
\sigma &=\partial \overline{f}+\overline{\partial }f \\
&=\sum_{j,\alpha }(f_j^{\overline{\alpha }}\theta ^j\otimes
\widetilde{\eta } _{\overline{\alpha }}+f_{\overline{j}}^\alpha
\theta ^{\overline{j}}\otimes \widetilde{\eta }_\alpha).
\endaligned \tag{3.12}
$$
Similarly we have
$$
\aligned
\tau &=\partial f+\overline{\partial f} \\
&=\sum_{j,\alpha }(f_j^\alpha \theta ^j\otimes \widetilde{\eta
}_\alpha +f_{ \overline{j}}^{\overline{\alpha }}\theta
^{\overline{j}}\otimes \widetilde{ \eta }_{\overline{\alpha }}).
\endaligned\tag{3.13}
$$
The covariant derivative of $\sigma $ is given by
$$
\aligned \nabla \sigma =&\sum_{j,k,\alpha
}\{f_{jk}^{\overline{\alpha }}\theta ^j\otimes \theta ^k\otimes
\widetilde{\eta }_{\overline{\alpha
}}+f_{j\overline{k}}^{\overline{\alpha }}\theta ^j\otimes \theta
^{\overline{k}}\otimes \widetilde{\eta }_{\overline{\alpha }} \\
&+f_{\overline{j}k}^\alpha \theta ^{\overline{j}}\otimes \theta
^k\otimes \widetilde{\eta }_\alpha
+f_{\overline{j}\overline{k}}^\alpha \theta ^{
\overline{j}}\otimes \theta ^{\overline{k}}\otimes \widetilde{\eta
}_\alpha \}.
\endaligned\tag{3.14}
$$
By (1.2) , we have
$$
\aligned
(d^\nabla \sigma )(X,Y)&=(\nabla _X\sigma )(Y)-(\nabla _Y\sigma )(X) \\
&=(\nabla \sigma )(Y;X)-(\nabla \sigma )(X;Y).
\endaligned\tag{3.15}
$$
By (3.14), (3.15) and using (2.31), we get
$$
\aligned d^\nabla \sigma =&\sum_{j,k,\alpha
}\{f_{jk}^{\overline{\alpha }}(\theta ^k\wedge \theta ^j)\otimes
\widetilde{\eta }_{\overline{\alpha }}+f_{j
\overline{k}}^{\overline{\alpha }}(\theta ^{\overline{k}}\wedge
\theta
^j)\otimes \widetilde{\eta }_{\overline{\alpha }} \\
&+f_{\overline{j}k}^\alpha (\theta ^k\wedge \theta
^{\overline{j}})\otimes \widetilde{\eta }_\alpha
+f_{\overline{j}\overline{k}}^\alpha (\theta ^{
\overline{k}}\wedge \theta ^{\overline{j}})\otimes \widetilde{\eta
}_\alpha
\} \\
=&\sum_{j,k,\alpha }(f_{k\overline{j}}^\alpha (\theta ^k\wedge
\theta ^{\overline{j}})\otimes \widetilde{\eta }_\alpha
-f_{j\overline{k}}^{\overline{\alpha }}(\theta ^j\wedge \theta
^{\overline{k}})\otimes \widetilde{\eta }_{\overline{\alpha }}) \\
=&\sum_{k,j,\alpha }(f_{k\overline{j}}^\alpha (\theta ^k\wedge
\theta ^{\overline{j}})\otimes \widetilde{\eta }_\alpha
-f_{k\overline{j}}^{\overline{\alpha }}(\theta ^k\wedge \theta
^{\overline{j}})\otimes \widetilde{\eta }_{ \overline{\alpha }}).
\endaligned\tag{3.16}
$$
By (2.8) and (2.9), we have $df=\sigma +\tau $. Since $d^\nabla
df=0$, it follows from (3.16) that
$$
d^\nabla \tau =-d^\nabla \sigma =\sum_{k,j,\alpha
}(-f_{k\overline{j} }^\alpha (\theta ^k\wedge \theta
^{\overline{j}})\otimes \widetilde{\eta } _\alpha
+f_{k\overline{j}}^{\overline{\alpha }}(\theta ^k\wedge \theta ^{
\overline{j}})\otimes \widetilde{\eta }_{\overline{\alpha }}).
\tag{3.17}
$$
From (3.12) and (3.16), we obtain
$$
\aligned & <(d^\nabla \sigma )(\eta _k,e_A),\sigma
(e_A)>\\
&=<(d^\nabla \sigma )(\eta _k,\eta _j),\sigma (\eta
_{\overline{j}})>
+<(d^\nabla \sigma )(\eta _k,\eta _{\overline{j}}),\sigma (\eta _j)> \\
&=\sum_{j,\alpha }f_{k\overline{j}}^\alpha f_j^{\overline{\alpha
}}\endaligned \tag{3.18}
$$
and
$$
<(d^\nabla \sigma )(\eta _{\overline{k}},e_A),\sigma
(e_A)>=\sum_{j,\alpha }f_{\overline{k}j}^{\overline{\alpha
}}f_{\overline{j}}^\alpha. \tag{3.19}
$$
Define a $1$-form $\gamma $ as follows
$$
\gamma (X)=\sum_A<(d^\nabla \sigma )(X,e_A),\sigma (e_A)>.
\tag{3.20}
$$
By (3.18) and (3.19), we obtain
$$
\aligned \gamma &=\sum_k[<(d^\nabla \sigma )(\eta _k,e_A),\sigma
(e_A)>\theta ^k+<(d^\nabla \sigma )(\eta
_{\overline{k}},e_A),\sigma (e_A)>\theta ^{
\overline{k}}] \\
&=\sum_k[(\sum_{j,\alpha }f_{k\overline{j}}^\alpha
f_j^{\overline{\alpha } })\theta ^k+\sum_{j,\alpha
}(f_{\overline{k}j}^{\overline{\alpha }}f_{ \overline{j}}^\alpha
)\theta ^{\overline{k}}].
\endaligned\tag{3.21}
$$
Consequently we get
$$
\aligned div(\gamma )&=\sum_k[(\sum_{j,\alpha
}f_{k\overline{j}}^\alpha f_j^{\overline{ \alpha
}})_{\overline{k}}+(\sum_{j,\alpha }(f_{\overline{k}j}^{\overline{
\alpha }}f_{\overline{j}}^\alpha ))_k] \\
&=2\sum_{k,j,\alpha }|f_{k\overline{j}}^\alpha
|^2+\sum_{k,j,\alpha }(f_{k \overline{j}\overline{k}}^\alpha
f_j^{\overline{\alpha }}+f_{\overline{k} jk}^{\overline{\alpha
}}f_{\overline{j}}^\alpha).
\endaligned\tag{3.22}
$$
By Ricci identity and K\"ahler identity, we have
$$
\aligned (f_{k\overline{j}\overline{k}}^\alpha
-f_{k\overline{k}\overline{j}}^\alpha ) \widetilde{\eta }_\alpha
&=-f_k^\beta \widetilde{R}_\beta ^\alpha
(f_{\overline{j}},f_{\overline{k}}) \widetilde{\eta }_\alpha
\endaligned\tag{3.23}
$$
where $\widetilde{R}$ denotes the curvature tensor of $N$. From
(2.29) and (3.23), it follows that
$$
\aligned f_{k\overline{j}\overline{k}}^\alpha
-f_{k\overline{k}\overline{j}}^\alpha &=-f_k^\beta
\widetilde{R}_{\beta \gamma \overline{\delta }}^\alpha f_{
\overline{j}}^\gamma f_{\overline{k}}^{\overline{\delta
}}-f_k^\beta \widetilde{R}_{\beta \overline{\gamma }\delta
}^\alpha f_{\overline{j}}^{
\overline{\gamma }}f_{\overline{k}}^\delta \\
f_{\overline{k}jk}^{\overline{\alpha
}}-f_{\overline{k}kj}^{\overline{\alpha
}}&=-f_{\overline{k}}^{\overline{\beta
}}\widetilde{R}_{\overline{\beta } \overline{\gamma }\delta
}^{\overline{\alpha }}f_j^{\overline{\gamma } }f_k^\delta
-f_{\overline{k}}^{\overline{\beta }}\widetilde{R}_{\overline{
\beta }\gamma \overline{\delta }}^{\overline{\alpha }}f_j^\gamma
f_k^{ \overline{\delta }}.
\endaligned
$$
Hence
$$
\aligned &div(\gamma
)\\=&2\sum_{k,j,\alpha}|f_{k\overline{j}}^\alpha
|^2-\sum_{k,j,\alpha }\widetilde{R}_{\beta \gamma \overline{\delta
}}^\alpha f_j^{\overline{\alpha }}f_k^\beta
f_{\overline{j}}^\gamma f_{\overline{k}}^{ \overline{\delta
}}-\sum_{k,j,\alpha }\widetilde{R}_{\beta \overline{\gamma }
\delta }^\alpha f_j^{\overline{\alpha }}f_k^\beta
f_{\overline{j}}^{\overline{\gamma }}f_{\overline{k}}^\delta \\
&-\sum_{k,j}\widetilde{R}_{\overline{\beta }\overline{\gamma
}\delta }^{ \overline{\alpha }}f_{\overline{j}}^\alpha
f_{\overline{k}}^{\overline{\beta }}f_j^{\overline{\gamma
}}f_k^\delta -\sum_{k,j}\widetilde{R}_{\overline{ \beta }\gamma
\overline{\delta }}^{\overline{\alpha }}f_{\overline{j} }^\alpha
f_{\overline{k}}^{\overline{\beta }}f_j^\gamma f_k^{\overline{
\delta }} \\
=&2\sum_{k,j,\alpha }|f_{k\overline{j}}^\alpha |^2-\sum_{k,j}\{\ll
\widetilde{Q}(f_j^{(0,1)}\wedge f_k^{(1,0)}),f_j^{(0,1)}\wedge f_k^{(1,0)}\gg \\
&+\ll \widetilde{Q}(f_j^{(0,1)}\wedge
f_k^{(1,0)}),f_j^{(1,0)}\wedge f_k^{(0,1)}\gg +\ll
\widetilde{Q}(f_j^{(0,1)}\wedge
f_k^{(1,0)}),f_j^{(0,1)}\wedge f_k^{(1,0)}\gg \\
&+\ll \widetilde{Q}(f_j^{(1,0)}\wedge
f_k^{(0,1)}),f_j^{(0,1)}\wedge
f_k^{(1,0)}\gg \} \\
=&2\sum_{k,j,\alpha }|f_{k\overline{j}}^\alpha
|^2-\sum_{k,j}\{2\ll \widetilde{Q}(f_k^{(1,0)}\wedge
f_j^{(0,1)}),f_k^{(1,0)}\wedge f_j^{(0,1)}\gg
\\
&-\ll \widetilde{Q}(f_k^{(1,0)}\wedge
f_j^{(0,1)}),f_j^{(1,0)}\wedge f_k^{(0,1)}\gg -\ll
\widetilde{Q}(f_j^{(1,0)}\wedge f_k^{(0,1)}),f_k^{(1,0)}\wedge
f_j^{(0,1)}\gg \}
\endaligned
$$
that is,
$$
div(\gamma )=2\sum_{k,j,\alpha }|f_{k\overline{j}}^\alpha
|^2-\sum_{k,j}\ll \widetilde{Q} (f_k\wedge f_j),f_k\wedge f_j\gg.
\tag{3.24}
$$
Likewise we define a $1$-form $\rho $ by
$$
\rho (X)=\sum_A<(d^\nabla \tau )(X,e_A),\tau (e_A)>.  \tag{3.25}
$$
From (3.13) and (3.17), we get
$$
\aligned
\rho &=<d^\nabla \tau (\eta _k,e_A),\tau (e_A)>\theta
^k+<d^\nabla \tau (\eta
_{\overline{k}},e_A),\tau (e_A)>\theta ^{\overline{k}} \\
&=\sum_k[(\sum_{j,\alpha }f_{k\overline{j}}^{\overline{\alpha
}}f_j^\alpha )\theta ^k+(\sum_{j,\alpha }(f_{\overline{k}j}^\alpha
f_{\overline{j}}^{\overline{\alpha }})\theta ^{\overline{k}}].
\endaligned\tag{3.26}
$$
A similar computation yields
$$
div(\rho )=2\sum_{k,j}|f_{k\overline{j}}^\alpha |^2-\sum_{k,j}\ll
\widetilde{ Q}(f_k\wedge f_j),f_k\wedge f_j\gg.\tag{3.27}
$$
Set $\widetilde{J^{\prime }}=-J^{\prime }$. Then we may define
$1-$forms $ \widetilde{\sigma }$ and $\widetilde{\tau }$ by (2.8),
(2.9) and using the new complex structure $\widetilde{J}^{\prime
}$ of $N$. It follows that $\sigma = \widetilde{\tau }$ and $\tau
=\widetilde{\sigma }$. The right hand side of (3.24) is obviously
independent of the choice of $J^{\prime }$ and $
\widetilde{J^{\prime }}$. This explains the result $div(\gamma
)=div(\rho )$ . Clearly either (3.24) or (3.27) can be used to
derive the pluriharmonicity and thus holomorphicity of a harmonic
map between two compact K\"ahler manifolds. In this paper, we try
to compute some terms related to the stress-energy tensors. These
two formulae, combined with (1.13), will become useful tools for
investigate harmonic maps between complete noncompact K\"ahler
manifolds too. From (3.24) and (3.27), we immediately have

\proclaim{Lemma 3.2} Let $f:M\rightarrow N$ be a harmonic map
between two K\"ahler manifolds. If $ N$ has strongly semi-negative
curvature, then $div(\gamma )=div(\rho )\geq 0$.
\endproclaim

\proclaim{Lemma 3.3} Let $M$ be a Riemannian manifold possing an
exhaustion function $\Phi $ with the properties (1.15) and (1.16).
Suppose $Z$ a vector field on $M$ such that
$$
\underset{R\rightarrow \infty }\to{\lim \inf}\frac 1R\int_{B_\Phi
(2R)-B_\Phi (R)}\Vert Z \Vert=0.\tag{3.28}
$$
If $divZ$ has an integral (that is, if either $(divZ)^{+}$ or
$(divZ)^{-}$ is integrable), then $\int_MdivZ=0$. In particular,
if outside some compact set $divZ$ is everywhere $\geq 0$ (or
$\leq 0$) then $\int_MdivZ=0$.
\endproclaim
\demo{Proof} This result was established in [Ka] for a complete
Riemannian manifold with $\Phi =r$ (the distance function). The
proof of Lemma 3.3 goes almost the same as the proof of the main
theorem in [Ka] , except for that one should use the general
exhaustion function $\Phi $ to replace the distance function. We
omit the detailed proof here. \qed
\enddemo
\remark{Remark 3.1} It is clear that if $\lim \inf_{R\rightarrow
\infty }\frac 1R\int_{B_\Phi (R)}||Z||=0$, then $Z$ satisfies
(3.28).
\endremark

\proclaim{Proposition 3.4} Let $M$ be a K\"ahler manifold possing
an exhaustion function $\Phi $ with the properties (1.15) and
(1.16). Let $f:M\rightarrow N$ be a harmonic map from a K\"ahler
manifold to a K\"ahler manifold with strongly semi-negative
curvature. If $f$ satisfies either
$$
\underset{R\rightarrow \infty}\to{\lim \inf} \frac 1R\int_{B_\Phi
(2R)-B_\Phi (R)}\Vert \gamma \Vert=0\tag{3.29}
$$
or
$$
\underset{R\rightarrow \infty}\to{\lim \inf} \frac 1R\int_{B_\Phi
(2R)-B_\Phi (R)}\Vert \rho \Vert=0\tag{3.30}
$$
then $f$ is pluriharmonic.
\endproclaim

\demo{Proof}Without loss of generality, we assume $f$ satisfies
(3.29). By Lemma 3.2, $div(\gamma )\geq 0$. Hence we obtain from
Lemma 3.3 and (3.24) the following
$$
\aligned \int_Mdiv(\gamma )&=2\sum_{k,j,\alpha
}\int_M|f_{k\overline{j}}^\alpha |^2-\sum_{k,j}\int_M\ll
\widetilde{Q}(f_k\wedge f_j),f_k\wedge f_j\gg \\
&=0.\endaligned\tag{3.31}
$$
This proposition follows immediately from (3.31).\qed
\enddemo

By (3.20) and (3.26), we get
$$
\aligned
\Vert\gamma \Vert &\leq |\overline{\partial }f|\cdot |(\nabla df)^{(1,1)}|\\
\Vert\rho \Vert &\leq |\partial f|\cdot |(\nabla df)^{(1,1)}|.
\endaligned
$$
Consequently we have

\proclaim{Corollary 3.5} Let $f:M\rightarrow N$ be a harmonic map
from a complete noncompact K\"ahler manifold to a K\"ahler
manifold with strongly semi-negative curvature. If $f$ satisfies
either
$$
\underset{r\rightarrow \infty }\to{\lim \inf }\frac
1r\int_{B_{2r}-B_r}|\overline{\partial }f|\cdot |(\nabla
df)^{(1,1)}|=0
$$
or
$$
\underset{r\rightarrow \infty }\to{\lim \inf }\frac
1r\int_{B_{2r}-B_r}|\partial f|\cdot |(\nabla df)^{(1,1)}|=0
$$
then $f$ is pluriharmonic.
\endproclaim

\proclaim{Theorem 3.6} Let $M$ be a K\"ahler manifold possing an
exhaustion function $\Phi $ with the properties (1.15) and (1.16).
Let $f:(M,g)\rightarrow (N,h)$ be a harmonic map into a K\"ahler
manifold of strongly semi-negative curvature. If $f$ satisfies
either
$$
(\int_{\partial B_\Phi (s)}|\overline{\partial }f|^2)^{-1}\notin
L^1(+\infty )\tag{3.32}
$$
or
$$
(\int_{\partial B_\Phi (s)}|\partial f|^2)^{-1}\notin L^1(+\infty
)\tag{3.33}
$$
then $f$ is pluriharmonic. In particular, $f$ satisfies (3.9),
that is,
$$
\ll \widetilde{Q}(f_j\wedge f_k),f_j\wedge f_k\gg =0,\quad
j,k=1,...,m \tag{3.34}
$$
where $\widetilde{Q}$ denotes the curvature operator of $N$.
\endproclaim

\demo{Proof} Without loss of generality, we assume that $f$
satisfies (3.32). It follows from (3.24) and the divergence
theorem that
$$
\aligned \int_{B_\Phi (t)}\sum_{k,j}|f_{k\overline{j}}^\alpha
|^2&\leq \frac
12\int_{B_\Phi (t)}div(\gamma ) \\
&=\frac 12\int_{\partial B_\Phi (t)}i_\nu \gamma
\endaligned\tag{3.35}
$$
where $\nu$ denotes the unit outward normal vector field of
$\partial B_\Phi (t)$. By (3.21), we have
$$
|i_\nu \gamma |\leq 2|\overline{\partial }f|(\sum_{k,j,\alpha }|f_{k%
\overline{j}}^\alpha |^2)^{1/2}
$$
which yields
$$
\int_{\partial B_\Phi (t)}i_\nu \gamma \leq 2\left\{
\int_{\partial B_\Phi (t)}|\overline{\partial }f|^2\right\}
^{1/2}\left\{ \int_{\partial B_\Phi (t)}\sum_{k,j,\alpha
}|f_{k\overline{j}}^\alpha |^2\right\} ^{1/2}\tag{3.36}
$$
Set
$$
\beta (t)=\int_{B_\Phi (t)}\sum_{k,j,\alpha
}|f_{k\overline{j}}^\alpha |^2
$$
Then, by the co-area formula,
$$
\aligned \beta ^{\prime }(t)&=\frac d{dt}\int_0^t(\int_{\partial
B_\Phi (s)}\frac{\sum_{k,j,\alpha }|f_{k\overline{j}}^\alpha |^2}{|\nabla \Phi |} )ds\\
&=\int_{\partial B_\Phi (t)}\frac{\sum_{k,j,\alpha
}|f_{k\overline{j}}^\alpha|^2}{|\nabla \Phi |} \\
&\geq \frac 1{\sqrt{\alpha }}\int_{\partial B_\Phi
(t)}\sum_{k,j,\alpha }|f_{k \overline{j}}^\alpha |^2
\endaligned\tag{3.37}
$$
where $\alpha $ is defined by (1.15). Putting (3.35) , (3.36) and
(3.37) together and squaring we finally get
$$
\beta (t)^2\leq \sqrt{\alpha }\left( \int_{\partial B_\Phi
(t)}|\overline{\partial }f|^2\right) \beta ^{\prime }(t)
\tag{3.38}
$$
Suppose that $(\nabla df)^{(1,1)}\neq 0$. Therefore there exists
$t_0>0 $ sufficiently large such that $\beta (t)>0$, for every
$t\geq t_0$. Fix such an $t_0$. From (3.38) we deduce the
following
$$
\sqrt{\alpha }(\beta (t_0)^{-1}-\beta (t)^{-1})\geq
\int_{t_0}^t\frac{ds}{\int_{\partial B_\Phi
(s)}|\overline{\partial }f|^2}
$$
and letting $t\rightarrow +\infty $ we contradict (3.32). Hence
$f$ is pluriharmonic. By definition, $\gamma \equiv 0$. Then
(3.24) implies that $f$ satisfies (3.9). \qed
\enddemo

\proclaim{Corollary 3.7} Let $M$ be a K\"ahler manifold possing an
exhaustion function $\Phi $ with the properties (1.15) and (1.16).
Let $f:(M,g)\rightarrow (N,h) $ be a harmonic map into a K\"ahler
manifold with strongly semi-negative curvature. If $f$ satisfies
either
$$
\int_{B_\Phi (R)}|\overline{\partial }f|^2\leq CR^2  \tag{3.39}
$$
or
$$
\int_{B_\Phi (R)}|\partial f|^2\leq CR^2  \tag{3.40}
$$
for any $R>0$ and some constant $C>0$, then $f$ is pluriharmonic.
In particular, $f$ satisfies (3.34).
\endproclaim
\demo{Proof} Set
$$
h(t)=\int_{B_\Phi (t)}|\overline{\partial }f|^2  \tag{3.41}
$$
So, by co-area formula,
$$
h^{\prime }(t)=\int_{\partial B_\Phi (t)}\frac{|\overline{\partial
}f|^2}{ |\nabla \Phi |}\geq \frac 1{\sqrt{\alpha }}\int_{\partial
B_\Phi (t)}| \overline{\partial }f|^2  \tag{3.42}
$$
Fix $t_0>0$, and let $t>t_0$. From Proposition 1.3 of [RS], we
know that
$$
\left. \int_{t_0}^t\left( \frac{s-t_0}{h(s)}\right) ds\leq
C_1\int_{t_0}^t \frac{ds}{h^{\prime }(s)}\right.   \tag{3.43}
$$
for some constant $C_1>0$. Then (3.42) and (3.43) imply
$$
\int_{t_0}^t\frac{s-t_0}{\int_{B_\Phi (s)}|\overline{\partial
}f|^2}ds\leq
C_1\sqrt{\alpha }\int_{t_0}^t\frac{ds}{\int_{\partial B_\Phi (s)}|\overline{%
\partial }f|^2}  \tag{3.44}
$$
Suppose $f$ satisfies (3.39). This implies
$$
\frac s{\int_{B_\Phi (s)}|\overline{\partial }f|^2}\notin
L^1(+\infty ) \tag{3.45}
$$
Thus we deduce from (3.44) and (3.45) that $f$ satisfies (3.32).
Likewise we may prove that (3.40) implies (3.33). Hence we prove
this corollary.\qed
\enddemo
\remark{Remark 3.2} Suppose $f:M\rightarrow N$ is a harmonic map
from a complete K\"ahler manifold to a K\"ahler manifold with
strongly semi-negative curvature. We may take the distance
function $r$ as the exhaustion function in Corollary 3.7.
Therefore we have proved that if one of the partial energies of
$f$ has growth order $\leq 2$ w.r.t. the distance function $r$,
then $f$ is pluriharmonic. In particular, if either $E^{\prime
}(f)<+\infty$ or $E^{\prime \prime }(f)<+\infty $, then $f$ is
pluriharmonic. In [Li] and [PRS], the authors gave similar
conditions about energy (not the partial energies) to derive the
pluriharmonicity of a harmonic map from a complete K\"ahler
manifold to a Riemannian manifold of non-positive Hermitian
curvature. Obviously any criteria for pluriharmonicity is
superfluous if $dim_CM=1$. It would be interesting to know whether
the growth order condition in Corollary 3.7 is optimal for
ensuring the pluriharmonicity or not.
\endremark

\heading{\bf 4. Monotonicity formulae of harmonic maps}
\endheading
\vskip 0.3true cm {\bf 4.1. The case of pluriharmonic maps} \vskip
0.3 true cm

First we show that the $1$-forms $\sigma $ and $\tau $ defined in
$\S 2$ satisfy the conservation laws if $f:M\rightarrow N$ is
pluriharmonic.

\proclaim{Lemma 4.1} A map $f:M\rightarrow N$ between two K\"ahler
manifolds is pluriharmonic if and only if the $1$-form $\sigma \in
A^1(f^{-1}TN)$ defined by (2.8) is closed or equivalently, the
$1$-form $\tau \in A^1(f^{-1}TN)$ defined by (2.9) is closed.
\endproclaim

\demo{Proof} By definition, we have
$$
\aligned &d^\nabla\sigma (X_1,X_2)\\&=\nabla _{X_1}\sigma
(X_2)-\nabla _{X_2}\sigma
(X_1)-\sigma ([X_1,X_2]) \\
&=\frac 12\{(\nabla df)(X_2,X_1)-(\nabla df)(X_1,X_2)+J^{\prime
}[(\nabla
df)(X_1,JX_2)-(\nabla df)(X_2,JX_1)]\} \\
&=\frac 12J^{\prime }[(\nabla df)(X_1,JX_2)-(\nabla df)(X_2,JX_1)]
\endaligned
$$
Likewise we may get
$$
d^\nabla \tau (X_1,X_2)=-\frac 12J^{\prime }[(\nabla
df)(X_1,JX_2)-(\nabla df)(X_2,JX_1)]
$$
It follows from Lemma 2.2 that $f$ is pluriharmonic iff $d^\nabla
\sigma =0$ iff $d^\nabla \tau =0$. \qed \enddemo

\proclaim{Proposition 4.2} Suppose $f:M\rightarrow N$ is a
pluriharmonic map between K\"ahler manifolds. Then $divS_\sigma
=divS_\tau =0$.
\endproclaim

\demo{Proof} This proposition follows immediately from Lemma 1.1,
Lemma 2.1, Lemma 4.1. \qed
\enddemo

\remark{Remark 4.1} Proposition 4.2 and (2.20) yield that $div\Psi
_f=0$ for a pluriharmonic map $f$, where $\Psi _f$ is the
$2-$tensor field defined by (2.19).
\endremark

From Proposition 1.2 and Proposition 4.2, we may deduce the energy
monotonicity formulae for $\sigma $ and $\tau $ of pluriharmonic
maps provided that $M$ posses suitable exhaustion functions. From
now on we assume $\dim _CM=m\geq 2$, unless otherwise indicated.

Let $\Phi $ be a Lipschitz continuous function on $M^m$, which
satisfies the following conditions (cf. also [Ta2]):

(4.1) $\Phi \geq 0$ and $\Phi $ is an exhaustion function of $M$,
i.e., each sublevel set $B_\Phi (t):=\{\Phi <t\}$ is relatively
compact in $M$ for $t\geq 0$;

(4.2) $\Psi =\Phi ^2$ is of class $C^\infty $ and $\Psi $ has only
discrete critical points;

(4.3) The constant $k_1=\inf_{x\in M}\sum_{i=1}^{m-1}\varepsilon
_i(x)$ is positive where $\varepsilon _1\leq \varepsilon _2\leq
\cdots \leq \varepsilon _m$ are the eigenvalues of the complex
Hessian $H(\Psi)=(\Psi _{i\overline{j}})$. The constant
$k_2=\sup_{x\in M}|\nabla \Phi |^2$ is finite. Set
$$
\lambda =k_1/k_2  \tag{4.4}
$$
The function $\Phi $ with the properties (4.1), (4.2) and (4.3)
will be called {\sl a special exhaustion function}. Notice that
(4.1) implies that $Im(\Phi )=[0,+\infty )$. For our purpose, we
only consider the unbounded exhaustion function for a complex
manifold in this paper, although not all results need this
assumption.  In addition, the condition (4.3) for $\Psi $ is
stronger than the condition for a function to be strictly
$(m-1)$-plurisubharmonic (see \S 6 for some discussion about the
notion of strict $(m-1)$-plurisubharmonicity).

\proclaim{Theorem 4.3} Let $f:M\rightarrow N$ be a pluriharmonic
map from a complete K\"ahler manifold to a K\"ahler manifold.
Suppose $M$ posses a special exhaustion function $\Phi $
satisfying (4.1), (4,2) and (4.3). Then
$$
\frac 1{\rho _1^\lambda }\int_{B_\Phi (\rho
_1)}|\overline{\partial } f|^2\leq \frac 1{\rho _2^\lambda
}\int_{B_\Phi (\rho _2)}|\overline{\partial }f|^2  \tag{4.5}
$$
and
$$
\frac 1{\rho _1^\lambda }\int_{B_\Phi (\rho _1)}|\partial f|^2\leq
\frac 1{\rho _2^\lambda }\int_{B_\Phi (\rho _2)}|\partial f|^2
\tag{4.6}
$$
for any $0<\rho _1\leq \rho _2$, where $\lambda $ is given by
(4.4).
\endproclaim

\demo{Proof} Take $X=\frac 12\nabla \Psi =\Phi \nabla \Phi $.
Obviously $(\nabla
\Psi )|_{\partial B_\Phi (t)}$ is an outward normal vector field along $%
\partial B_\Phi (t)$ for a regular value $t>0$ of $\Phi $. Thus $\nabla \Psi
=w(x)\nu $ on $\partial B_\Phi (t)$ with $w(x)>0$ for each point
$x\in
\partial B_\Phi (t)$, where $\nu $ denotes the unit outward normal vector
field of $\partial B_\Phi (t)$. By the definition of $S_\sigma $,
we have
$$
\aligned
 S_\sigma (X,\nu )&=\frac{|\sigma |^2}2<X,\nu >_g-<\sigma
(X),\sigma (\nu )>_h
\\
&=t\frac{|\sigma |^2}2<\nabla \Phi ,\nu >_g-\frac w2|\sigma (v)|_h^2 \\
&\leq t\sqrt{k_2}|\overline{\partial }f|^2
\endaligned\tag{4.7}
$$
on $\partial B_\Phi (t)$ and
$$
\aligned
<S_\sigma ,\nabla \theta _X>&=\frac 12<S_\sigma ,Hess(\Psi )> \\
&=\frac 12[\frac{|\sigma |^2}2\bigtriangleup _g\Psi -<\sigma \odot
\sigma ,Hess(\Psi )>]
\endaligned \tag{4.8}
$$
where $\bigtriangleup _g$ denotes the Laplace-Beltrami operator on
$M$. We choose a unitary basis $\{\eta
_i=(e_i-iJe_i)/\sqrt{2}\}_{i=1,...,m}$ at a point $p\in B_\Phi
(t)$ such that
$$
Hess(\Psi )(\eta _i,\overline{\eta }_j(p))=\varepsilon _i\delta
_{ij}
$$
which is equivalent to
$$
\aligned
Hess(\Psi )(e_i,e_j)+Hess(\Psi )(Je_i,Je_j)&=2\varepsilon _i\delta _{ij} \\
Hess(\Psi )(e_i,Je_j)-Hess(\Psi )(Je_i,e_j)&=0
\endaligned\tag{4.9}
$$
Obviously $\{e_i,Je_i\}_{i=1,...,m}$ is an orthonormal basis. Then
(4.9) gives
$$
\bigtriangleup _g\Psi =2\sum_{i=1}^m\varepsilon _i  \tag{4.10}
$$
Using (2.12) and (4.9), we obtain
$$
\aligned &<\sigma \odot \sigma ,Hess(\Psi )>\\
&=\sum_{i,j}\{<\sigma (e_i),\sigma(e_j)>Hess(\Psi )(e_i,e_j)+<\sigma (Je_i),\sigma (Je_j)>Hess(\Psi )(Je_i,Je_j) \\
&\quad+<\sigma (e_i),\sigma (Je_j)>Hess(\Psi )(e_i,Je_j)+<\sigma
(Je_i),\sigma(e_j)>Hess(\Psi )(Je_i,e_j)\} \\
&=\sum_{i,j}<\sigma (e_i),\sigma (e_j)>[Hess(\Psi
)(e_i,e_j)+Hess(\Psi)(Je_i,Je_j)] \\
&\quad+\sum_{i,j}<\sigma (e_i),\sigma (Je_j)>[Hess(\Psi
)(e_i,Je_j)-Hess(\Psi)(Je_i,e_j)] \\
&=2\sum_{i=1}^m|\sigma (e_i)|^2\varepsilon _i
\endaligned\tag{4.11}
$$
From (4.8), (4.10) and (4.11), we get
$$
\aligned <S_\sigma ,\nabla \theta _X>&=\frac 12\{\frac{|\sigma
|^2}2(2\sum_{j=1}^m
\varepsilon _j)-2\sum_{i=1}^m|\sigma (e_i)|^2\varepsilon _i\} \\
&\geq \frac 12\{|\sigma |^2\sum_{j=1}^m\varepsilon
_j-2\sum_{i=1}^m|\sigma
(e_i)|^2\varepsilon _m\} \\
&\geq \frac{k_1|\sigma |^2}2=k_1|\overline{\partial }f|^2
\endaligned\tag{4.12}
$$
It follows from (1.14), (4.7) and (4.12) that
$$
\sqrt{k_2}t\int_{\partial B_\Phi (t)}|\overline{\partial }f|^2\geq
k_1\int_{B_\Phi (t)}|\overline{\partial }f|^2  \tag{4.13}
$$
The remaining arguments are similar to those in the proof of
Proposition 1.2. Using the coarea formula yields, we get
$$
\frac d{dt}\int_{B_\Phi (t)}|\overline{\partial }f|^2\geq \frac 1{
\sqrt{k_2}}\int_{\partial B_\Phi (t)}|\overline{\partial
}f|^2\tag{4.14}
$$
and thus
$$
\frac{\frac d{dt}\int_{B_\Phi (t)}|\overline{\partial
}f|^2}{\int_{B_\Phi (t)}|\overline{\partial }f|^2}\geq \frac
\lambda t
$$
for any $t>0$, where $\lambda =k_1/k_2$. By integration over
$[\rho _1,\rho _2]$, we have
$$
\frac 1{\rho _1^\lambda }\int_{B_\Phi (\rho
_1)}|\overline{\partial } f|^2\leq \frac 1{\rho _2^\lambda
}\int_{B_\Phi (\rho _2)}|\overline{\partial }f|^2
$$
for any $0<\rho _1\leq \rho _2$. Similarly we can prove
$$
\frac 1{\rho _1^\lambda }\int_{B_\Phi (\rho _1)}|\partial f|^2\leq
\frac 1{\rho _2^\lambda }\int_{B_\Phi (\rho _2)}|\partial f|^2
$$
for any $0<\rho _1\leq \rho _2$. \qed
\enddemo

\remark{Remark 4.2} \newline (a) Suppose $\widehat{\lambda }$ is a
positive number less than $\lambda $ , that is,
$0<\widehat{\lambda }$ $<\lambda $. Clearly we have the
corresponding monotonicity formulae by replacing $\lambda $ by
$\widehat{\lambda }$ in (4.5) and (4.6). The larger the growth
order $\lambda$ we get, the better the monotonicity formulae
become.
\newline (b) In [Ta2], K. Takegoshi derived similar monotonicty
formulae for the energy of pluriharmonic maps. Here we establish
monotonicity formulae for the partial energies.
\endremark

Now we give some examples of K\"ahler manifolds which poss the
special exhaustion functions (Some of them were also discussed in
[Ta2] with somewhat different notations) .

\example{Example 4.1}\newline (a) Let $C^m$ be an $m\geq 2$
complex Euclidean space with the canonical K\"ahler metric. Take
$\Psi =\sum_{j=1}^mz_jz_{\overline{j}}=||z||^2$ and $ \Phi
=||z||$. Then $(\Psi _{i\overline{j}})=2(\delta _{ij})_{m\times
m}$ where the complex Hessian is diagonalized w.r.t. $\{\eta
_j=\frac 1{\sqrt{2} }(e_j-iJe_j)\}_{j=1}^m$. By definition, we
have $k_1=2m-2$ , $k_2=1$ and thus $\lambda =2m-2$.
\newline (b) Let $i:M^m\rightarrow C^N$
be an $m-$dimensional closed complex submanifold and
$F:C^N\rightarrow R$ be a smooth function on $C^N$. By the
composition formula of maps (cf.[EL]), we have
$$
Hess(F\circ i)(X,Y)=(HessF)(X,Y)+dF(B(X,Y))
$$
for any $X,Y\in TM$, which yields
$$
Hess(F\circ i)(\eta _i,\overline{\eta }_j)=(HessF)(\eta
_i,\overline{\eta } _j)+dF(B(\eta _i,\overline{\eta }_j))
$$
where $\{\eta _i=\frac 1{\sqrt{2}}(e_i-iJe_i)\}_{i=1}^m$ is any
unitary frame tangent to $M$. Since $i(M)$ is a complex
submanifold, we have $B(\eta _i,\overline{\eta }_j)=0$. Then
$$
Hess(F\circ i)(\eta _j,\overline{\eta }_k)=(HessF)(\eta
_j,\overline{\eta } _k)
$$
Now let $F=||z||^2=\sum_{A=1}^N|z_A|^2$ and set $\Psi =\Phi
^2=i^{*}(F)=F\circ i$. Obviously $k_1=2(m-1)$ and $k_2\leq 1$. If
necessary, translating the original point to a general position,
we may assume that $F$ is a Morse function (cf. [Mi]). So $F$ has
only discrete critical points. Then $\Phi $ is a special
exhaustion function with $\lambda \geq 2m-2$. Recall that every
Stein manifold $M^m$ can be realized as a closed submanifold of
$C^N$ by a proper holomorphic map $\psi :M^m\rightarrow C^N$. Thus
a Stein manifold $M$ admits a special exhaustion $\Psi =\Phi ^2$
with $\Phi =\psi ^{*}(||z||)$ and $\lambda \geq 2m-2$. It is known
that every closed complex submanifold of a Stein manifold is a
Stein manifold too. Therefore Stein manifolds provide us many
examples of K\"ahler manifolds which poss special exhaustion
functions.
\endexample

Notice that the special exhaustion functions in (b) of example 4.1
are obtained from extrinsic distance functions. Next we will show
that under suitable curvature conditions, the distance functions
of K\"ahler manifolds are special exhaustion functions too.

\proclaim{Lemma 4.4} Let $M^m$ be an $m$ dimensional complete
K\"ahler manifold with a pole $x_0$. Let $\varepsilon _1\leq
\varepsilon _2\leq \cdots \leq \varepsilon _m$ be the eigenvalues
of the complex Hessian $((r^2)_{i\overline{j}})$, where $r$ is the
distance function relative to $x_0$. Suppose there exists a
positive function $h(r)$ on $(0,+\infty )$ such that
$$
Hess(r)\geq h(r)[g-dr\otimes dr]
$$
then
$$
\sum_{i=1}^{m-1}\varepsilon _i\geq \cases
1+(2m-3)rh(r) & \text{if}\quad rh(r)\geq 1 \\
2(m-1)rh(r) & \text{if}\quad rh(r)<1
\endcases
$$
\endproclaim
\demo{Proof} By assumptions, we have
$$
Hess(r^2)\geq 2dr\otimes dr+2rh(r)[g-dr\otimes dr] \tag{4.15}
$$
First we consider the case $rh(r)\geq 1$. Replacing $\Psi $ by
$r^2$ in (4.9), we see that $2$ $\sum_{i=1}^{m-1}\varepsilon _i$
is the trace of $ Hess(r^2)$ on some real $2(m-1)$ dimensional
$J-$invariant subspace of $TM$. It follows from (4.15) that
$$
\sum_{i=1}^{m-1}\varepsilon _i\geq 1+(2m-3)rh(r)
$$
Clearly we have
$$
\sum_{i=1}^{m-1}\varepsilon _i\geq 2(m-1)rh(r)
$$
if $rh(r)<1$. This proves this lemma. \qed
\enddemo

\proclaim{Lemma 4.5} Let $(M,g)$ be a complete Riemannian manifold
with a pole $x_0$ and $r$ the distance function relative to $x_0$.
Denote by $K_r$ the radial curvature $K_r$ of $M$.

(i) If $K_r\leq 0$, then
$$
Hess(r)\geq \frac 1r[g-dr\otimes dr]
$$

(ii) If $K_r\leq \frac{b^2}{1+r^2}$ with $b^2\in [0,1/4]$, then
$$
Hess(r)\geq \frac{1+\sqrt{1-4b^2}}{2r}[g-dr\otimes dr]
$$

(iii) If $K_r\leq \frac B{(1+r^2)^{1+\varepsilon }}$ with
$\varepsilon >0$ and $0\leq B<2\varepsilon $, then
$$
Hess(r)\geq \frac{1-\frac B{2\varepsilon }}r[g-dr\otimes dr]
$$

(iv) If $K_r\leq -\beta ^2$ with $\beta >0$, then
$$
Hess(r)\geq \beta \coth (\beta r)[g-dr\otimes dr]
$$

(v) If $K_r\leq -\frac{a^2}{1+r^2}$ with $a>0$, then
$$
Hess(r)\geq \max \{\frac{1+\sqrt{1+4a^2}}{2(1+r)},\frac
1r\}[g-dr\otimes dr]
$$
\endproclaim

\demo{Proof} The cases (i) and (iv) are standard (cf. [GW]). The
case (ii) is proved in Lemma 1.2 (b) of [EF]. The case (iii)
follows immediately from the quasi-isometry Theorem due to [GW]
(cf. also [DW]). The case (v) is treated in [GW], [PRS] as an
asymptotical comparison theorem. Actually we may deduce (cf. page
39 of [PRS])
$$
Hess(r)\geq \frac{1+\sqrt{1+4a^2}}{2(1+r)}[g-dr\otimes dr]
$$
On the other hand, the assumption $K_r\leq -\frac{a^2}{1+r^2}<0$
implies by (i) that
$$
Hess(r)\geq \frac 1r[g-dr\otimes dr]
$$
Therefore we prove (v). \qed
\enddemo

\proclaim{Lemma 4.6}
Let $M^m$ be an $m$ dimensional complete K\"ahler manifold with a pole $x_0$%
. Let $\varepsilon _1\leq \varepsilon _2\leq \cdots \leq
\varepsilon _m$ be the eigenvalues of $((r^2)_{i\overline{j}})$,
where $r$ is the distance function relative to $x_0$. Suppose the
radial curvature $K_r$ of $M$ satisfies one of the five conditions
(i), (ii), (iii), (iv) and (v) in Lemma 4.5. Then
$$
\sum_{i=1}^{m-1}\varepsilon _i\geq \cases 2m-2 & \text{if}\quad
K_r\text{ satisfies (i)}\\
(m-1)(1+\sqrt{1-4b^2}) & \text{if}\quad K_r\text{ satisfies (ii)}\\
2(m-1)(1-\frac B{2\varepsilon }) & \text{if}\quad K_r\text{
satisfies (iii)}\\
1+(2m-3)\beta r\coth (\beta r) & \text{if}\quad K_r\text{ satisfies (iv)}\\
\max \{2m-2,1+\frac{(2m-3)(1+\sqrt{1+4a^2})r}{2(1+r)}\} &
\text{if}\quad K_r\text{  satisfies (v)}\endcases \tag{4.16}
$$
\endproclaim

\demo{Proof} Suppose $K_r$ satisfies (i). Then Lemma 4.5 and Lemma
4.4 yield immediately $\sum_{i=1}^{m-1}\varepsilon _i\geq  2m-2$.
This proves (i). Now we assume that $K_r$ satisfies (iv). Clearly
$\beta r\coth (\beta r)\geq 1$ on $(0,+\infty )$, because the
increasing function $\beta r\coth \beta r\rightarrow 1$ as
$r\rightarrow 0$. Therefore Lemma 4.4 implies that
$$
\sum_{i=1}^{m-1}\varepsilon _i\geq 1+(2m-3)\beta r\coth (\beta r)
$$
Hence we have proved (iv). In a similar way, we may use Lemma 4.4
and Lemma 4.5 to prove the cases (ii), (iii) and (v) too. \qed
\enddemo

From Lemma 4.6, we see that if $M$ is as in Lemma 4.6, then its
distance function $r$ is a special exhaustion function for $M$.
From the proof of Theorem 4.3 and remark 4.2, it follows that

\proclaim{Theorem 4.7}
Let $M$ ($\dim _CM=m>1$), $r$, $K_r$ be as in Lemma 4.6. Suppose $%
f:M\rightarrow N$ is a pluriharmonic map between K\"ahler
manifolds. Set
$$
\lambda =\cases
2m-2 & \text{if}\quad K_r\text{ satisfies one of (i), (iv), (v)} \\
(m-1)(1+\sqrt{1-4b^2}) & \text{if}\quad K_r\text{ satisfies (ii)} \\
2(m-1)(1-\frac B{2\varepsilon }) & \text{if}\quad K_r\text{
satisfies (iii)}
\endcases \tag{4.17}
$$
Then $f$ satisfies
$$
\frac 1{\rho _1^\lambda }\int_{B_{\rho
_1}(x_0)}|\overline{\partial } f|^2\leq \frac 1{\rho _2^\lambda
}\int_{B_{\rho _2}(x_0)}|\overline{\partial }f|^2
$$
and
$$
\frac 1{\rho _1^\lambda }\int_{B_{\rho _1}(x_0)}|\partial
f|^2\leq \frac 1{\rho _2^\lambda }\int_{B_{\rho _2}(x_0)}|\partial
f|^2
$$
for any $0<\rho _1\leq \rho _2$.
\endproclaim

Now let $f:M\rightarrow N$ be a pluriconformal harmonic map from a
K\"ahler manifold. By definition, $f$ satisfies
$$
<df(JX),df(JY)>=<df(X),df(Y)>\tag{4.18}
$$
for any $X,Y\in TM$. Clearly (4.18) is equivalent to
$$
<df(X),df(JY)>=-<df(JX),df(Y)>
$$

\proclaim{Theorem 4.8} Let $f:M\rightarrow N$ be a pluriconformal
harmonic map between two K\"ahler manifolds. Suppose $M$ satisfies
the conditions in Theorem 4.3 (resp. Theorem 4.7) and $\lambda $
is given by (4.4) (resp. (4.17)). Then $f$ satisfies
$$
\frac 1{\rho _1^\lambda }\int_{B_\Phi (\rho _1)}|df|^2\leq \frac
1{\rho _2^\lambda }\int_{B_\Phi (\rho _2)}|df|^2 \quad\text{(resp.
}\frac 1{\rho _1^\lambda }\int_{B\rho _1}|df|^2\leq \frac 1{\rho
_2^\lambda }\int_{B_{\rho _2}}|df|^2\text{)}
$$
for any $0<\rho _1\leq \rho _2$.
\endproclaim
\demo{Proof} Suppose $M$ posses a special function $\Phi $. As in
the proof of Theorem 4.3, we set $X=\Phi \nabla \Phi $. Since $f$
is harmonic, it is known that $f$ satisfies a conservation law
([BE]), that is, $divS_f=0$, where $S_f$ is defined by (1.6). Then
(1.13) yields
$$
\int_{\partial B_\Phi (t)}S_f(X,\nu )=\int_{B_\Phi (t)}<S_f,\nabla
\theta _X>
$$
for a regular value $t>0$ of $\Phi $. Similar to (4.7), we get
$$
S_f(X,\nu )\leq \frac{t\sqrt{k_2}}2|df|^2
$$
By using the pluriconformality of $f$ and replacing $\sigma $ by
$df$ in (4.8), (4.11) and (4.12), we may derive the following:
$$
<S_f,\nabla \theta _X>\geq \frac{k_1}2|df|^2
$$
It follows that
$$
t\sqrt{k_2}\int_{\partial B_\Phi (t)}|df|^2\geq k_1\int_{B_\Phi
(t)}|df|^2
$$
Similar to the remaining argument in the proof of Theorem 4.3, we
have
$$
\frac 1{\rho _1^\lambda }\int_{B_\Phi (\rho _1)}|df|^2\leq \frac
1{\rho _2^\lambda }\int_{B_\Phi (\rho _2)}|df|^2
$$
for any $0<\rho _1\leq \rho _2$.\qed
\enddemo

Notice that the curvature tensors of the target manifolds in
Theorem 4.3, Theorem 4.7 and Theorem 4.8 play no role, since we
consider the conservative cases in this subsection.

\vskip 0.3true cm {\bf 4.2. The case of harmonic maps} \vskip 0.3
true cm

Now we consider a harmonic map between K\"ahler manifolds. In this
case, $\sigma $ and $\tau $ don't satisfy the conservation laws in
general. However, if the target K\"ahler manifold has strongly
semi-negative curvature, we will show that the integral formula
(1.13) can still be used to establish the monotonicity formulae
for the partial energies. \proclaim{Lemma 4.9} Let $f:M\rightarrow
N$ be a harmonic map from a K\"ahler manifold into a K\"ahler
manifold with strongly semi-negative curvature. Let $D\subset M$
be a domain with a compact closure and non-empty smooth boundary.
Let $\varphi $ be a smooth defining function for $D$ with only
discrete critical point. Set $X=\nabla \varphi $. Then
$$
\int_D(divS_\sigma )(X)=\int_D(divS_\tau )(X)\geq 0  \tag{4.19}
$$
\endproclaim

\demo{Proof} From Lemma 1.1 and Lemma 2.1, we have
$$
(divS_\sigma )(X)=<i_Xd^\nabla \sigma ,\sigma >  \tag{4.20}
$$
Set $D_t=\{\varphi <t\}$ for $t\leq 0$. Let $\nu =\frac{\nabla
\varphi }{ |\nabla \varphi |}$ on $\partial D_t$, where $t$ is a
regular value in $ Im(\varphi )$. By (4.20) and using the coarea
formula and divergence theorem, we deduce that
$$
\aligned \int_D(divS_\sigma )(X)&=\int_{-\infty }^0(\int_{\partial
D_t}(divS_\sigma )(\frac{\nabla \varphi }{|\nabla \varphi |}) \\
&=\int_{-\infty }^0(\int_{\partial D_t}<i_\nu d^\nabla \sigma ,\sigma >)dt \\
&=\int_{-\infty }^0(\int_{\partial D_t}i_\nu \gamma )dt \\
&=\int_{-\infty }^0(\int_{D_t}div\gamma )dt
\endaligned
$$
where $\gamma $ is defined by (3.20). Likewise we have
$$
\int_D(divS_\tau )(X)=\int_{-\infty }^0(\int_{D_t}div\rho )dt
$$
where $\rho $ is defined by (3.25). It follows from Lemma 3.2 that
$ \int_D(divS_\sigma )(X)=\int_D(divS_\tau )(X)\geq 0$.\qed
\enddemo

\remark{Remark 4.3}\newline (a) From the proof of Lemma 4.9, we
see that the corresponding results also hold true if the integral
domain $D$ in (4.19) is replaced by the domain $ D(c,d):=\{x\in
D:c<\varphi <d\}$ for any $c<d\leq 0$.
\newline
(b) If $M$ is a complete K\"ahler manifold with a pole $x_0$, we
may take $\varphi =\frac 12r^2$ , where $r$ is distance function
relative $x_0$ . Set $ X=r\frac \partial {\partial r}$ in Lemma
4.9. Then we have
$$
\int_{B_R}(divS_\sigma )(X)=\int_{B_R}(divS_\tau )(X)\geq 0
\tag{4.21}
$$
\newline
(c) From (3.24) and(3.27), we see that if $f$ is harmonic map, but
not pluriharmonic. Then $div(\gamma )=div(\rho )>0$ at some point
of the geodesic ball $B_R$. Therefore $\int_{B_R}div(S_\sigma
)(X)=\int_{B_R}(divS_\tau )(X)>0$, which implies that $\sigma $
and $\tau $ don't satisfy the conservation laws.
\endremark

\proclaim{Theorem 4.10} Let $f:M\rightarrow N$ be a harmonic map
from a complete K\"ahler manifold to a K\"ahler manifold with
strongly semi-negative curvature. Suppose $M$ satisfies the
conditions in Theorem 4.3 (resp. Theorem 4.7) , $\lambda $ is
given by (4.4) (resp. (4.17)) and set $D_R=B_\Phi (R)$ (resp.
$B_R(x_0)$). Then
$$
\frac 1{\rho _1^\lambda }\int_{D_{\rho _1}}|\overline{\partial
}f|^2\leq \frac 1{\rho _2^\lambda }\int_{D_{\rho
_2}}|\overline{\partial }f|^2 \tag{4.22}
$$
and
$$
\frac 1{\rho _1^\lambda }\int_{D_{\rho _1}}|\partial f|^2\leq
\frac 1{\rho _2^\lambda }\int_{D_{\rho _2}}|\partial f|^2
\tag{4.23}
$$
for any $0<\rho _1\leq \rho _2$.
\endproclaim

\demo{Proof} Without loss of generality, we assume that $M$
satisfies the conditions in Theorem 4.3. Actually we have already
pointed out this result in Remark 1.1. Here we give only a brief
discussion. Take $X=\nabla (\frac 12\Psi )=\Phi \nabla \Phi $.
Since $N$ has strongly semi-negative curvature, we obtain by Lemma
4.9 that
$$
\int_{B_\Phi (R)}divS_\sigma (X)\geq 0  \tag{4.24}
$$
It follows from (1.13) and (4.24) that
$$
\int_{\partial B_\Phi (R)}S_\sigma (X,\nu )\geq \int_{B_\Phi
(R)}<S_\sigma ,\nabla \theta _X>
$$
The remaining argument is similar to that in the proof of Theorem
4.3. \qed
\enddemo

\vskip 3true cm\vskip 0.3true cm {\bf 4.3. Local monotonicity
formulae} \vskip 0.3 true cm

\proclaim{Theorem 4.11} Let $M$ be an $m-$dimensional K\"ahler
manifold and let $r$ denote the distance function relative to
$x_0\in M$. Suppose the radial curvature $K_r$ of $M$ satisfies
$K_r\leq K_0$ on $B_{R_0}(x_0)$, where $K_0$ is a positive
constant and $R_0$ is a fixed positive number less than the
injective radius of $M$ at $x_0$. Suppose $f:M\rightarrow N$ is
either a pluriharmonic map into any K\"ahler manifold or a
harmonic map into a K\"ahler manifold with strongly semi-negative
curvature. Then
$$
\frac{e^{C(K_0)\rho _2}}{\rho _2^{2m-2}}\int_{B_{\rho
_2}(x_0)}|\overline{\partial }f|^2-\frac{e^{C(K_0)\rho _1}}{\rho
_1^{2m-2}}\int_{B_{\rho _1}(x_0)}| \overline{\partial }f|^2\geq
\int_{\rho _1}^{\rho _2}[\frac{e^{C(K_0)t}}{
t^{2m-2}}\int_{\partial B_t(x_0)}|i_{\frac \partial {\partial
r}}\sigma |^2]dt
$$
and
$$
\frac{e^{C(K_0)\rho _2}}{\rho _2^{2m-2}}\int_{B_{\rho
_2}(x_0)}|\partial f|^2- \frac{e^{C(K_0)\rho _1}}{\rho
_1^{2m-2}}\int_{B_{\rho _1}(x_0)}|\partial f|^2\geq \int_{\rho
_1}^{\rho _2}[\frac{e^{C(K_0)t}}{t^{2m-2}}\int_{\partial
B_t(x_0)}|i_{\frac
\partial {\partial r}}\tau |^2]dt
$$
for any $0<\rho _1\leq \rho _2<R_0$, where $C(K_0)$ is a constant
depending on $K_0$.
\endproclaim

\demo{Proof} Since $K\leq K_0$ on $B_{R_0}(p)$, we get by Hessian
comparison theorem the following
$$
Hess(r)\geq \sqrt{K_0}\cot (\sqrt{K_0}r)[g-dr\otimes dr]
$$
Now the strictly decreasing function $\sqrt{K_0}r\cot (\sqrt{K_0}
r)\rightarrow 1$ as $r\rightarrow 0$, we have $\sqrt{K_0}r\cot
(\sqrt{K_0} r)<1$ on $(0,+\infty )$. Let $\Psi =r^2$ and let
$\varepsilon _1\leq \cdots \leq \varepsilon _m$ be the eigenvalue
of $H(r^2)$. It follows from Lemma 4.4 that
$$
\aligned
\sum_{i=1}^{m-1}\varepsilon _i&\geq 2(m-1)\sqrt{K_0}r\cot (\sqrt{K_0}r) \\
&=2(m-1)+2(m-1)[\sqrt{K_0}r\cot (\sqrt{K_0}r)-1]
\endaligned\tag{4.25}
$$
Obviously there exists a constant positive $\widetilde{C}(K_0)$
such that
$$
1-\sqrt{K_0}r\cot (\sqrt{K_0}r)\leq r\widetilde{C}(K_0) \tag{4.26}
$$
on $B_{R_0}(p)$. Set $X=\frac 12\nabla \Psi $. From (4.12), (4.25)
and (4.26), we get
$$
\aligned <S_\sigma ,\nabla \theta _X>&\geq
(\sum_{j=1}^{m-1}\varepsilon _j)\frac{|\sigma |^2}2 \\
&\geq [2(m-1)-C(K_0)r]\frac{|\sigma |^2}2
\endaligned\tag{4.27}
$$
where $C(K_0)=2(m-1)\widetilde{C}(K_0)$.

Notice that $|\nabla r|=1$. Then (1.8) gives
$$
S_\sigma (X,\nu )\leq r[\frac{|\sigma |^2}2-|i_\nu \sigma |^2]
\tag{4.28}
$$
When $f$ is pluriharmonic, we know from Proposition 4.2 that
$\sigma $ satisfy the conservation law. Using a similar technique
as in the proof of Proposition 1.2, we can deduce from (4.27) and
(4.28) the following
$$
\frac d{dr}[e^{C(K_0)r}r^{-(2m-2)}\int_{B_r(x_0)}|\sigma |^2]\geq
2r^{-(2m-2)}e^{C(K_0)r}\int_{\partial B_r(x_0)}|\sigma |^2
\tag{4.29}
$$
for $r<R_0$. By integration on $[\rho _1,\rho _2]$, we get the
monotonicity formula for $\sigma $. Likewise we have the
monotonicity formula for $\tau $.

Suppose now that $f:M\rightarrow N$ is a harmonic map into a
K\"ahler manifold with strongly semi-negative curvature. From
Lemma 4.9, Remark 4.3 and Remark 1.1, it is clear that the
monotonicity formulae still hold. \qed
\enddemo

\proclaim{Corollary 4.12} Let $B_{R_0}(x_0)\subset M$, $N$ and
$f:M\rightarrow N$ be as in Theorem 4.11. Then
$$
\frac{e^{C(\alpha ,\beta )\rho _1}}{\rho _1^{2m-2}}\int_{B_{\rho
_1}(x_0)}| \overline{\partial }f|^2\leq \frac{e^{C(\alpha ,\beta
)\rho _2}}{\rho _2^{2m-2}}\int_{B_{\rho
_2}(x_0)}|\overline{\partial }f|^2
$$
and
$$
\frac{e^{C(\alpha ,\beta )\rho _1}}{\rho _1^{2m-2}}\int_{B_{\rho
_1}(x_0)}|\partial f|^2\leq \frac{e^{C(\alpha ,\beta )\rho
_2}}{\rho _2^{2m-2}} \int_{B_{\rho _2}(x_0)}|\partial f|^2
$$
for $0<\rho _1\leq \rho _2<R_0$.
\endproclaim

Next, we hope to establish the monotonicity formulae outside of a
compact subset of a complete K\"ahler manifold.

\proclaim{Theorem 4.13} Let $M$ be an $m-$dimensional complete
K\"ahler manifold. Suppose $M$ posses an exhaustion function $\Phi
$ which is special outside a sublevel set $ B_\Phi (R_0)$ for some
$R_0>0$. Suppose $f:M\rightarrow N$ is either a pluriharmonic map
into any K\"ahler manifold or a harmonic map into a K\"ahler
manifold with strongly semi-negative curvature. Then
$$
\frac 1{\rho _1^\lambda }\int_{B_\Phi (\rho _1)-B_\Phi
(R_0)}|\overline{\partial }f|^2\leq \frac 1{\rho _2^\lambda
}\int_{B_\Phi (\rho _2)-B_\Phi (R_0)}|\overline{\partial }f|^2
\tag{4.30}
$$
and
$$
\frac 1{\rho _1^\lambda }\int_{B_\Phi (\rho _1)-B_\Phi
(R_0)}|\partial f|^2\leq \frac 1{\rho _2^\lambda }\int_{B_\Phi
(\rho _2)-B_\Phi (R_0)}|\partial f|^2  \tag{4.31}
$$
for any $R_0<\rho _1\leq \rho _2$, where $\lambda $ is defined by
(4.4) on $M-B_\Phi (R_0)$.
\endproclaim

\demo{Proof} Take $X=\frac 12\nabla \Phi ^2=\Phi \nabla \Phi $.
For any $R>R_0$, we set $D=B_\Phi (R)-B_\Phi (R_0)$. By applying
the integral formula (1.13) on $D$ and arguing in a similar way as
in Proposition 1.2, Theorem 4.3 and Theorem 4.10, we may deduce
the following
$$
\aligned &\int_{\partial B_\Phi (R)}S_\sigma (X,\nu
)-\int_{\partial B_\Phi
(R_0)}S_\sigma (X,\nu ) \\
&=R<\nabla \Phi ,\nu >\{\int_{\partial B_\Phi (R)}S_\sigma (\nu
,\nu
)-\int_{\partial B_\Phi (R_0)}S_\sigma (\nu ,\nu )\} \\
&\geq \frac 12\int_{B_\Phi (R)-B_\Phi (R_0)}<S_\sigma ,Hess(\Phi ^2)> \\
&\geq \frac{k_1}2\int_{B_\Phi (R)-B_\Phi (R_0)}|\sigma |^2
\endaligned\tag{4.32}
$$
By definition of $S_\sigma $ and (2.12), we get
$$
\aligned
S_\sigma (\nu ,\nu )&=\frac{|\sigma |^2}2-<\sigma (\nu ),\sigma (\nu )> \\
&\geq \frac 12[<\sigma (\nu ),\sigma (\nu )>+<\sigma (J\nu
),\sigma (J\nu )>]\\
&-<\sigma (\nu ),\sigma (\nu )> \\
&=0
\endaligned\tag{4.33}
$$
Then (4.32) and (4.33) yield
$$
R\sqrt{k_2}\int_{\partial B_\Phi (R)}S_\sigma (\nu ,\nu )\geq
\frac{k_1} 2\int_{B_\Phi (R)-B_\Phi (R_0)}|\sigma |^2  \tag{4.34}
$$
It follows from the coarea formula and (4.34) that
$$
\aligned R\frac d{dR}(\int_{B_\Phi (R)-B_\Phi (R_0)}|\sigma
|^2)&\geq \frac R{\sqrt{k_2}}\int_{\partial B_\Phi (R)}|\sigma |^2 \\
&\geq \frac{2R}{\sqrt{k_2}}\int_{\partial B_\Phi (R)}S_\sigma (\nu ,\nu ) \\
&\geq \lambda \int_{B_\Phi (R)-B_\Phi (R_0)}|\sigma |^2
\endaligned
$$
which implies that
$$
\frac d{dr}\{R^{-\lambda }\int_{B_\Phi (R)-B_\Phi
(R_0)}|\overline{\partial } f|^2\}\geq 0  \tag{4.35}
$$
By integrating (4.35) on $[\rho _1,\rho _2]$, we get (4.30).
Likewise we have (4.31). \qed
\enddemo

As an application of Theorem 4.20, we give the following

\proclaim{Theorem 4.14} Let $M$ be a complete K\"ahler manifold
with a pole $x_0$. Suppose the radial curvature $K_r$ of $M$
satisfies one of the following two conditions

(a) $K_r$ $\leq -\beta ^2$ with $\beta >0$;

(b) $K_r\leq -\frac{a^2}{1+r^2}$ with $a>0$.
\newline
Suppose $f:M\rightarrow N$ is either a pluriharmonic map into a
K\"ahler manifold or a harmonic map into a K\"ahler manifold with
strongly semi-negative curvature. Set
$$
\lambda _0=\cases
1+(2m-3)\beta R_0\coth (\beta R_0) &\text{ if }K_r\text{ satisfies (a)} \\
C(R_0) &\text{ if }K_r\text{ satisfies (b)}
\endcases
$$
for any $R_0>0$, where
$$
C(R_0)=\cases 1+\frac{(2m-3)(1+\sqrt{1+4a^2})R_0}{2(1+R_0)} &
\text{if }\frac{(1+\sqrt{1+4a^2})R_0}{2(1+R_0)}\geq 1 \\
\frac{(m-1)(1+\sqrt{1+4a^2})R_0}{(1+R_0)} & \text{if
}\frac{(1+\sqrt{1+4a^2} )R_0}{2(1+R_0)}<1
\endcases
$$
Then
$$
\frac 1{\rho _1^{\lambda _0}}\int_{B_{\rho _1}(x_0)-B_{R_0}(x_0)}|\overline{%
\partial }f|^2\leq \frac 1{\rho _2^{\lambda _0}}\int_{B\rho
_2(x_0)-B_{R_0}(x_0)}|\overline{\partial }f|^2  \tag{4.36}
$$
and
$$
\frac 1{\rho _1^{\lambda _0}}\int_{B_{\rho
_1}(x_0)-B_{R_0}(x_0)}|\partial f|^2\leq \frac 1{\rho _2^{\lambda
_0}}\int_{B\rho _2(x_0)-B_{R_0}(x_0)}|\partial f|^2  \tag{4.37}
$$
for any $R_0<\rho _1\leq \rho _2$.
\endproclaim

\demo{Proof} For the case $K_r\leq -\beta ^2$, we have from Lemma
4.4, Lemma 4.5 and the proof of Lemma 4.6, we have
$\sum_{i=1}^{m-1}\varepsilon _i(x)\geq 1+(2m-3)\beta R_0\coth
(\beta R_0)$ for $x\in M-B_{R_0}(x_0)$. Then Theorem 4.13 follows
immediately from Theorem 4.14. In a similar way, we may prove the
case (b) by Lemma 4.4 and Lemma 4.5 and using the fact that the
function $\frac r{1+r}$ is increasing. \qed
\enddemo

Notice that (4.33) is the key property that allows us to establish
the monotonicity formulae outside a sublevel set of the exhaustion
function $\Phi$ in Theorem 4.13. Suppose now that $f:M\rightarrow
N$ is a pluriconformal harmonic map from a K\"ahler manifold. We
may consider the stress energy tensor $S_f$ as in Theorem 4.8. Let
$\nu $ be the unit outward normal vector field of $B_\Phi (R_0)$.
By (4.18), we also have
$$
S_f(\nu ,\nu )\geq 0  \tag{4.38}
$$
Similar to the arguments in Theorem 4.13 (see also Theorem 4.8),
it is easy to deduce from (4.38) the following:

\proclaim{Theorem 4.15} Let $M$, $\Phi $, $R_0$ and $\lambda $ be
as in Theorem 4.13. Suppose $f:M\rightarrow N$ is a pluriconformal
harmonic map. Then
$$
\frac 1{\rho _1^\lambda }\int_{B_\Phi (\rho _1)-B_\Phi
(R_0)}|df|^2\leq \frac 1{\rho _2^\lambda }\int_{B_\Phi (\rho
_2)-B_\Phi (R_0)}|df|^2
$$
for any $R_0<\rho _1\leq \rho _2$.
\endproclaim

\proclaim{Corollary 4.16} Let $M$, $K_r$, $R_0$ and $\lambda _0$
be as in Theorem 4.14. Suppose $f:M\rightarrow N$ is a
pluriconformal harmonic map. Then
$$
\frac 1{\rho _1^{\lambda _0}}\int_{B_{\rho _1}-B_{R_0}}|df|^2\leq
\frac 1{\rho _2^{\lambda _0}}\int_{B_{\rho _2}-B_{R_0}}|df|^2
$$
for any $R_0<\rho _1\leq \rho _2$.
\endproclaim

\heading{\bf 5. Holomorphicity and constancy of harmonic maps}
\endheading
\vskip 0,3 true cm

In this section, we derive some results about holomorphicity and
constancy of harmonic maps between K\"ahler manifolds. Most of
these results are direct consequences of the monotonicity formulae
in last section.

\proclaim{Theorem 5.1}Let $M$, $\Phi $, $R_0$ and $\lambda $ be as
in Theorem 4.13, that is, $\Phi $ is an exhaustion function on
$M$, which is special on $M-B_\Phi (R_0)$. Suppose $f:M\rightarrow
N$ is either a pluriharmonic map into a K\"ahler manifold or a
harmonic map into a K\"ahler manifold with strongly semi-negative
curvature. If
$$\int_{B_\Phi (R)}|\overline{\partial
}f|^2=o(R^\lambda )\quad\text{(resp. }\int_{B_\Phi (R)}|\partial
f|^2=o(R^\lambda ))\quad \text{ as }R\rightarrow \infty
\tag{5.1}$$ then $f$ is holomorphic (resp. anti-holomorphic). In
particular, if $E^{\prime \prime }(f)<+\infty $ (resp. $E^{\prime
}(f)<+\infty $), then $f$ is holomorphic (resp. anti-holomorphic).
\endproclaim
\demo{Proof} The assumption $\int_{B_\Phi (R)}|\overline{\partial
}f|^2=o(R^\lambda )$ (resp. $\int_{B_\Phi (R)}|\partial
f|^2=o(R^\lambda )$) as $R\rightarrow \infty $ is equivalent to
$$
\int_{B_\Phi (R)-B_\Phi (R_0)}|\overline{\partial
}f|^2=o(R^\lambda )\quad\text{ (resp.}\int_{B_\Phi (R)}|\partial
f|^2=o(R^\lambda )\text{) as }R\rightarrow \infty
$$
By (4.30) (resp.(4.31)), we deduce that $f$ is holomorphic (resp.
anti-holomorphic) on $M-B_\Phi (R_0)$. Since $f$ is harmonic, it
then follows as in [Si1] from the unique continuation property
that $\overline{\partial }f=0$ (resp. $\partial f=0$) on the whole
$M$. This proves the theorem. \qed
\enddemo

\proclaim{Corollary 5.2} Let $M$, $\Phi $, $\lambda $ be as in
Theorem 4.3. Suppose $f:M\rightarrow N$ is either a pluriharmonic
map into a K\"ahler manifold or a harmonic map into a K\"ahler
manifold with strongly semi-negative curvature. If $f$ satisfies
the growth condition as in Theorem 5.1, then $f$ is holomorphic
(resp. anti-holomorphic).
\endproclaim

\proclaim{Corollary 5.3} Let $M$, $\Phi $, $\lambda $, $f$ and $N$
be as in Corollary 5.2. If $f$ satisfies
$$
\int_{B_\Phi (R)}|df|^2=o(R^\lambda )\quad\text{as}\quad
R\rightarrow \infty\tag{5.2}
$$
then $f$ is constant. In particular, if $E(f)<+\infty $, then $f$
is constant.
\endproclaim

\demo{Proof} By Corollary 5.2, $f$ is both holomorphic and
anti-holomorphic. Hence $f$ must be constant.\qed \enddemo

According to [Wu], a function on a manifold is called {\sl
quasipositive} if it is everywhere nonnegative and is positive at
one point.

\proclaim{Proposition 5.4} Let $M$ be a K\"ahler manifold possing
a function $\Phi $ with the properties (4.1), (4.2) and that
$\sum_{i=1}^{m-1}\varepsilon _i$ is quasipositive, where
$\varepsilon _1\leq \varepsilon _2\leq \cdots \leq \varepsilon _m$
are the eigenvalues of the complex Hessian $H(\Phi ^2)$. Suppose
$f:M\rightarrow N$ is either a pluriharmonic map into a K\"ahler
manifold or a harmonic map into a K\"ahler manifold with strongly
semi-negative curvature. If
$$
\int_M|\overline{\partial }f|^2<\infty \quad\text{ ( resp.
}\int_M|\partial f|^2<\infty \text{)}  \tag{5.3}
$$
then $f$ is holomorphic (resp. anti-holomorphic). In particular,
if $\int_M|df|^2<\infty$ then $f$ is constant.
\endproclaim

\demo{Proof} First we assume that $f$ is a pluriharmonic map
satisfying $E^{\prime \prime }(f)<\infty $. Set $X=\frac 12\nabla
(\Phi ^2)$. Similar arguments as in the proof of Theorem 4.3 yield
$$
<S_\sigma ,\nabla \theta _X>\geq (\sum_{i=1}^{m-1}\varepsilon
_i)|\overline{\partial }f|^2  \tag{5.4}
$$
From Proposition 4.2, (1.12) and (5.4), we obtain
$$
div(i_XS_\sigma )\geq (\sum_{i=1}^{m-1}\varepsilon
_i)|\overline{\partial } f|^2  \tag{5.5}
$$
Clearly $|i_XS_\sigma |\leq C(k_2)\Phi |\overline{\partial }f|^2$
for some constant $C(k_2)$. Thus
$$
\frac 1R\int_{B_\Phi (2R)-B_\Phi (R)}|i_XS_\sigma |\leq
2C(k_2)\int_{B_\Phi (2R)-B_\Phi (R)}|\overline{\partial }f|^2
$$
By $E^{\prime \prime }(f)<\infty $, we get
$$
\underset{R\rightarrow \infty }\to{\lim \inf \frac 1R}\int_{B_\Phi
(2R)-B_\Phi (R)}|i_XS_\sigma |=0
$$
Then Lemma 3.3 yields
$$
\int_Mdiv(i_XS_\sigma )=0  \tag{5.6}
$$
It follows from (5.5) and (5.6) that
$$
\int_M(\sum_{i=1}^{m-1}\varepsilon _i)|\overline{\partial }f|^2=0
$$
Since $\sum_{i=1}^{m-1}\varepsilon _i$ is quasipositive,
$|\overline{\partial }f|^2$ vanishes on an open subset of $M$.
Consequently $\overline{\partial }f=0$ on the whole $M$, that is,
$f$ is holomorphic.

Suppose now that  $f:M\rightarrow N$ is a harmonic map into a
K\"ahler manifold. If $E^{\prime \prime }(f)<\infty $, then
Corollary 3.7 implies that $f$ is pluriharmonic. The remaining
arguments are obvious. Thus we have proved this proposition. \qed
\enddemo

\proclaim{Theorem 5.5} Let $M$ be a complete K\"ahler manifold
with a pole $x_0$. Suppose the radial curvature $K_r$ of $M$
satisfies one of the following three conditions:

(i) $K_r\leq 0$;

(ii) $K_r\leq \frac{b^2}{1+r^2}$ with $b\in (0,1/4]$;

(iii) $K_r\leq \frac B{(1+r^2)^{1+\varepsilon }}$ with
$0<B<2\varepsilon $.
\newline Suppose $f:M\rightarrow N$ is either
a pluriharmonic map into a K\"ahler manifold or a harmonic map
into a K\"ahler manifold with strongly semi-negative curvature. If
$$
\int_{B_r}|\overline{\partial }f|^2=o(r^\lambda )\quad\text{
(resp. } \int_{B_r}|\partial f|^2=o(r^\lambda))\quad\text{as
}r\rightarrow \infty\tag{5.7}
$$
where
$$
\lambda =\cases
2m-2 & \text{if }K_r\text{ satisfies (i)} \\
(m-1)(1+\sqrt{1-4b)} & \text{if }K_r\text{ satisfies (ii)} \\
2(m-1)(1-\frac B{2\varepsilon }) & \text{if }K_r\text{ satisfies
(iii)}
\endcases\tag{5.8}
$$
then $f$ is holomorphic (resp. anti-holomorphic). In particular,
if $E^{\prime \prime }(f)<+\infty $ (resp. $E^{\prime }(f)<+\infty
$), then $f$ is holomorphic (resp. anti-holomorphic).
\endproclaim

\demo{Proof}This follows directly from Theorem 4.7 and Theorem
4.10. \qed \enddemo

\remark{Remark 5.1} The author in [Wa] proved that any harmonic
map $f:C^m\rightarrow C^n$ with $E^{\prime \prime }(f)<+\infty $
is holomorphic.
\endremark

\proclaim{Theorem 5.6} Let $M$ be a complete K\"ahler manifold
whose radial curvature $K_r$ satisfies one of the following two
conditions

(iv) $K_r$ $\leq -\beta ^2$ with $\beta >0$;

(v) $K_r\leq -\frac{a^2}{1+r^2}$ with $a>0$.
\newline
Suppose $f:M\rightarrow N$ is either a pluriharmonic map into a
K\"ahler manifold or a harmonic map into a K\"ahler manifold with
strongly semi-negative curvature. If
$$
\int_{B_r(x_0)}|\overline{\partial }f|^2=o(r^{\Lambda
_0})\quad\text{as }r\rightarrow \infty \tag{5.9}
$$
where $\Lambda _0$ is either any positive number if $K_r$
satisfies (iv) or any positive number less than
$1+\frac{(2m-3)(1+\sqrt{1+4a^2})}2$ if $K_r$ satisfies (v), then
$f$ is holomorphic. In particular, if $ E^{\prime \prime
}(f)<+\infty $, then $f$ is holomorphic. Likewise, if we replace
$\overline{\partial }f$ by $\partial f$ in the previous
conditions, then $f$ is anti-holomorphic.
\endproclaim

\demo{Proof} First we assume that $K_r$ satisfies (iv) and $f$
satisfies (5.9). Since $x\coth x$ is a nondecreasing function for
$x\geq 0$ and $\lim_{x\rightarrow +\infty }\coth x=1$, there
exists an $R_0$ such that $1+(2m-3)\beta R_0\coth (\beta R_0)\geq
\Lambda _0$. Then Theorem 4.14 implies that $f$ is holomorphic.

Next we assume that $K_r$ satisfies (v) and $f$ satisfies (5.9)
too. Notice that $\frac r{1+r}$ is an increasing function for
$r\geq 0$ and $ \lim_{r\rightarrow \infty }\frac r{1+r}=1$.
Clearly $\frac{(1+\sqrt{1+4a^2} )R_0}{2(1+R_0)}>1$ and
$1+\frac{(2m-3)(1+\sqrt{1+4a^2})R_0}{2(1+R_0)}\geq \Lambda _0$ for
a sufficiently large $R_0$. Hence we get from Theorem 4.14 that
$f$ is holomorphic. \qed
\enddemo

\remark{Remark 5.2} Let $f:M\rightarrow N$ be a harmonic map from
any complete K\"ahler manifold into a K\"ahler manifold with
strongly semi-negative curvature. Notice that we assert in
Corollary 3.7 (see also Remark 3.2), without assuming any
curvature conditions on $M$ , that if one of the partial energy
has growth order less than or equal to $2 $ then $f$ is
pluriharmonic. Suppose now that the radial curvature $K_r(M)$
satisfies one of the conditions in Theorem 4.7. We assert above
that if the partial energy has growth order less than $\lambda $,
then $f$ is $\pm $holomorphic. In most cases, $\lambda $ given by
(4.17) may be larger than $2$. However, it is easy to verify that
$\lambda $ is less than or equal to $2$ in following cases:
$$
\cases
m=2 & \text{if }K_r\text{ satisfies one of (i), (ii)}\\
2\leq m\leq 1+\frac 1{1-\frac B{2\varepsilon }} & \text{if
}K_r\text{ satisfies (iii)}
\endcases
$$
Remarkably we may conclude for the case (iv) of Theorem 5.6 that
if the $E''-$energy (resp. $E'-$energy) has polynomial growth in
$r$, then $f$ is holomorphic (resp. antiholomorphic).
\endremark

\proclaim{Corollary 5.7} Let $M$, $N$, $f$, $\lambda $ and
$\Lambda _0$ be as in Theorem 5.5 and Theorem 5.6. Suppose
$f:M\rightarrow N$ is either a pluriharmonic map into a K\"ahller
manifold or a harmonic map into a K\"ahler manifold with strongly
semi-negative curvature. If $f$ satisfies
$$
\int_{B_r}|df|^2=\cases
o(r^\lambda ) & \text{if }M\text{ is as in Theorem5.5} \\
o(r^{\Lambda _0}) & \text{if }M\text{ is as in Theorem5.6}
\endcases\tag{5.10}
$$
Then $f$ is constant. In particular, if $E(f)<+\infty $, then $f$
is constant.
\endproclaim
\demo{Proof} It follows immediately from Theorem 5.5 and Theorem
5.6 that $f$ is both holomorphic and anti-holomorphic. Hence $f$
must be constant.\qed
\enddemo

\proclaim{Theorem 5.8}
Let $M$ be as in Theorem 4.15 (resp. Corollary 4.16). Suppose $%
f:M\rightarrow N$ is a pluriconformal harmonic map from the
K\"ahler manifold $M$. If
$$
\int_{B_\Phi (R)}|df|^2=o(R^\lambda )\quad\text{ (resp.
}\int_{B_r}|df|^2=o(r^{ \lambda _0})\text{)}
$$
then $f$ is constant. In particular, if $E(f)<+\infty $, then $f$
is constant.
\endproclaim

\demo{Proof} This theorem follows immediately from Theorem 4.15
(resp. Corollary 4.16) and the unique continuation theorem of
harmonic maps (cf. [EL]). \qed
\enddemo

\remark{Remark 5.3} \newline (a) For most of the results in this
paper, we have to assumed $m>1$ for the conditions on $M$. We
cannot expect in general that the corresponding results hold in
the case of $m=1$. For example, there exit many harmonic maps from
$S^2\rightarrow CP^N$ which are neither holomorphic nor
anti-holomorphic (cf. [EW]). Since $R^2$ is conformally equivalent
to $S^2\backslash \{p\}$, it follows that there exit many harmonic
maps $R^2\rightarrow CP^N$ of finite energy, which are neither
holomorphic nor anti-holomorphic.
\newline (b) To establish Liouville Theorems for a harmonic
map $f$, one may also consider the stress-energy tensor $S_f$.
Usually some curvature pinching conditions on the domain manifolds
are needed to obtain the Liouville theorems under the energy
growth conditions (cf. [DW] for details). Note that in Corollary
5.7 or Theorem 5.8, we only assume some suitable upper bounds for
$K_r$ to establish the Liouville theorems.
\endremark

In Proposition 3.4, Corollary 3.5, Theorem 3.6 and Corollary 3.7,
we have deduced the pluriharmonicity of harmonic maps from any
complete K\"ahler manifold $M$ to a K\"ahler manifold with
strongly semi-negative curvature under some mild growth conditions
about the harmonic maps. In particular, all these harmonic maps
satisfy (3.34). When $N$ is an irreducible Hermitian symmetric
space, we define an integer $P(N)$ as follows (cf. [Si2]):
$$
P(N)=\cases
(p-1)(q-1)+1 & \text{ if }\widetilde{N}=D^{Ipq} \\
\frac 12(p-2)(p-3)+1 & \text{ if }\widetilde{N}=D^{IIp} \\
\frac 12p(p-1)+1 & \text{ if }\widetilde{N}=D^{IIIp} \\
2 & \text{ if }\widetilde{N}=D^{IVp} \\
6 & \text{ if }\widetilde{N}=D^V \\
11 & \text{ if }\widetilde{N}=D^{VI}
\endcases  \tag{5.11}
$$
where $\widetilde{N}$ denotes the universal covering space of $N$.
We have the following:

\proclaim{Lemma 5.9} (cf. [Si1,2]) Let $f:M\rightarrow N$ be a
smooth map between two K\"ahler manifolds. Then $f$ is holomorphic
or anti-holomorphic, provided (i) $N$ has strongly negative
curvature tensor, $f$ is a harmonic map with (3.34) and
$max_Mrank_Rdf\geq 4$, or (ii) $N$ is an irreducible Hermitian
symmetric space of noncompact type, $f$ is a pluriharmonic map
with $max_Mrank_Rdf\geq 2P(N)+1$, where $P(N)$ is given by
$(5.11)$.
\endproclaim

From Lemma 5.9, we obtain that:

\proclaim{Theorem 5.10} Let $f:M\rightarrow N$ be a harmonic map
from a complete K\"ahler manifold into a K\"ahler manifold.
Suppose $M$ and $f$ satisfy the conditions as in one of
Proposition 3.4, Corollary 3.5, Theorem 3.6 and Corollary 3.7.
Then $f$ is holomorphic or anti-holomorphic, provided that (i) $N$
has strongly negative curvature tensor and $max_Mrank_Rdf\geq 4$,
or (ii) $N$ is an irreducible Hermitian symmetric space of
noncompact type and $max_Mrank_Rdf\geq 2P(N)+1$.
\endproclaim

\heading{\bf 6. Harmonic maps with CR Dirichlet boundary-values}
\endheading
\vskip 0.3 true cm

Let $D$ $\subset M^m$ be a relatively compact domain with smooth
connected boundary $\partial D$. For $p\in \partial D$, we denote
by $H_p(\partial D)$ the real $ 2m-2$ dimensional subspace of
$T_p(\partial D)$ which is $J$ invariant. The distribution
$H:=\{H_p :  p\in \partial D\}$ on $\partial D$ is called the
holomorphic distribution of $\partial D$. Suppose $f:\partial
D\rightarrow N$ is a map into a K\"ahler manifold $N$. We say that
$f$ satisfies the tangential Cauchy-Riemann equation
$\overline{\partial }_bf=0$ on $\partial D $ if for every point
$p\in \partial D$, $\Pi ^{1,0}\circ df$ $(\xi )=0$ for any $\xi
\in H_p^C\cap T_p^{0,1}M$, where $\Pi ^{1,0}:TN\otimes
C\rightarrow T^{1,0}N$ is the natural projection. It is easy to
verify that $f$ satisfies the tangential Cauchy-Riemann equation
if and only if $\sigma $ annihilates any tangent vector in the
holomorphic distribution $H$.

\proclaim{Theorem 6.1} Let $\overline{D}_i\subset M$ ($i=1,2$) be
two connected, compact smooth domains in a K\"ahler manifold such
that $D_1\subset \subset D_2$. Set $ D=D_2-\overline{D}_1$.
Suppose there exists a function $\Psi \in C^2( \overline{D})$
satisfying the following properties:

(i) $k_1(x)=\sum_{i=1}^{m-1}\varepsilon _i(x)$ is quasipositive on
$D$, where $\varepsilon _1\leq \varepsilon _2\leq \cdots \leq
\varepsilon _m$ are the eigenvalues of the complex Hessian $(\Psi
_{i\overline{j}})$;

(ii) $\nabla \Psi =w(x)\nu \,$ on $\partial D_1$ with $w(x)\geq 0$
for every $x\in \partial D_1$, where $\nu $ is the unit outward
normal vector field along $\partial D_1$.
\newline
Suppose $f:D\rightarrow N$ is a pluriharmonic map into a K\"ahler
manifold $N$ such that $f\in C^2(\overline{D},N)$ and
$\overline{\partial } _bf=0$ on $\partial D_2$. Then $f$ is
holomorphic.
\endproclaim

\demo{Proof} Take $X=\frac 12\nabla \Psi $. Since $f:\overline{D}$
$\rightarrow N$ is a pluriharmonic map, it follows from
Proposition 4.2 that $f$ satisfies the conservation law. By
(1.14), we get
$$
\int_{\partial D_2}S_\sigma (X,\nu )-\int_{\partial D_1}S_\sigma
(X,\nu )=\int_D<S_\sigma ,\nabla \theta _X>  \tag{6.1}
$$
Similar to (4.12), we deduce that
$$
\aligned
<S_\sigma ,\nabla \theta _X>&=\frac 12<S_\sigma ,Hess(\Psi )> \\
&\geq k_1(x)|\overline{\partial }f|^2
\endaligned\tag{6.2}
$$
on $D$. For any $x\in \partial D_2$, we have
$$
\aligned 2S_\sigma (X,\nu )&=\frac{|\sigma |_x^2}2<\nabla \Psi
,\nu>_x-<\sigma (\nabla
\Psi ),\sigma (\nu )>_x \\
&=\frac{|\sigma |_x^2}2<\nabla \Psi ,\nu >_x-<\nabla \Psi ,\nu
>_x|\sigma (v)|_x^2-<\sigma((\nabla \Psi )_T),\sigma (\nu )>_x
\endaligned\tag{6.3}
$$
where $(\nabla \Psi )_T$ denotes the tangential projection of
$\nabla \Psi $ on $T(\partial D)$. Since $\overline{\partial
}_bf=0$ on $\partial D_2$, we have
$$
\sigma ((\nabla \Psi )_T)=<\nabla \Psi ,J\nu >\sigma (J\nu )
$$
From (2.12), we derive that
$$
|\sigma (\nu )|^2=|\sigma (J\nu )|^2  \tag{6.4}
$$
and
$$
\aligned  <\sigma (J\nu ),\sigma (\nu )>&=<\sigma (J^2\nu ),\sigma (J\nu )> \\
&=-<\sigma (\nu ),\sigma (J\nu )>\endaligned
$$
that is,
$$
<\sigma (J\nu ),\sigma (\nu )>=0  \tag{6.5}
$$
Therefore
$$
\aligned 2S_\sigma (X,\nu )&=<\nabla \Psi ,\nu >[\frac{|\sigma
|_q^2}2-|\sigma (\nu
)|^2] \\
&=<\nabla \Psi ,\nu >[\frac{|\sigma (\nu )|^2+|\sigma (J\nu
)|^2}2-|\sigma
(\nu )|^2] \\
&=0
\endaligned\tag{6.6}
$$
on $\partial D_2$. Similar to the argument in (4.33), we have the
following
$$
\aligned
S_\sigma (X,\nu )&=\frac w2S_\sigma (\nu ,\nu ) \\
&\geq 0
\endaligned\tag{6.7}
$$
on $\partial D_1$. From (6.1), (6.2), (6.6) and (6.7), we have
$$
\int_Dk_1(x)|\overline{\partial }f|^2\leq 0  \tag{6.8}
$$
By assumption $k_1(x)>0$ at some point $x\in D$, then (6.8)
implies that $\overline{\partial }f\equiv 0$ in a neighborhood $U$
of $p$. It follows that $f$ is holomorphic on $D$. \qed
\enddemo

It is clear that if $\overline{D}_1$ is a sublevel set of the
function $\Psi $, then the condition (ii) in Theorem 6.1 is
automatically satisfied.

\proclaim{Corollary 6.2} Let $M$ , $D_i$ ($i=1,2$) and $f$ be as
in Theorem 6.1. Suppose there exists a function $\Psi \in
C^2(\overline{D})$ satisfying the condition (i) of Theorem 6.1.
Suppose there is a real number $c$ such that $D_1\subset \{q\in
D:\Psi (q)\leq c\}\subset \subset D$. Then $f$ is holomorphic on
$D$.
\endproclaim
\demo{Proof} Set $\widehat{D}_1=\{q\in D:\Psi (q)<c\}$ and
$\widehat{D}=D_2- \widehat{D}_1$. Applying Theorem 6.1 to the
pluriharmonic map $f:\widehat{D} \rightarrow N$, we deduce that
$f$ is holomorphic on $\widehat{D}$. Hence $f$ is holomorphic on
$D$. \qed
\enddemo

Note also that if $D_1=\emptyset $, the condition (ii) of Theorem
6.1 is void. Therefore we get

\proclaim{Corollary 6.3} Let $\overline{D}\subset M$ be a compact
domain in a K\"ahler manifold with smooth connected boundary.
Suppose there exists a function $\Psi \in C^2( \overline{D})$
satisfying the property that the function
$k_1(x)=\sum_{i=1}^{m-1}\varepsilon _i(x)$ is quasipositive on
$D$, where $ \varepsilon _1\leq \varepsilon _2\leq \cdots \leq
\varepsilon _m$ are the eigenvalues of the complex Hessian $(\Psi
_{i\overline{j}})$. Suppose $f:D\rightarrow N$ is a pluriharmonic
map into a K\"ahler manifold $N$ such that $f\in
C^2(\overline{D},N)$ and $\overline{\partial } _bf=0$ on $\partial
D$. Then $f$ is holomorphic.
\endproclaim

\proclaim{Corollary 6.4} Let $M$ be as in Lemma 4.6 and let
$\overline{D}\subset M$ be any compact domain in $M$ with smooth
connected boundary. Suppose $f:D\rightarrow N$ is a pluriharmonic
map into a K\"ahler manifold $N$ such that $f\in
C^2(\overline{D},N)$ and $\overline{\partial } _bf=0$ on $\partial
D$. Then $f$ is holomorphic.
\endproclaim

\demo{Proof} Take $\Psi =r^2$. Then Lemma 4.6 implies that $\Psi $
satisfies the assumptions of Corollary 6.3 on any connected,
compact smooth domain $D$. Hence we prove this corollary. \qed
\enddemo
\remark{Remark 6.1} From the proof of Theorem 6.1, it is easy to
deduce the following result. Suppose $f:D\rightarrow N$ is a
pluriconformal harmonic map from $D$. If $df|_H=0$ for the
holomorphic distribution $H$ of $\partial D_2$, then $f$ is
constant. Here we have to use the unique continuation theorem of
harmonic maps. Hence we have the corresponding results of
Corollary 6.2, Corollary 6.3 and Corollary 6.4 for pluriconformal
harmonic maps. However, the conclusion for $f$ is the constancy
instead of the holomorphicity.
\endremark.

Notice that we don't assume any convexity conditions about
$\partial D$ and any curvature conditions on the target manifold
in Theorem 6.1, Corollary 6.2,  Corollary 6.3 and Corollary 6.4.
In addition, the conditions about $ \Psi $ in Corollary 6.3 are
weaker than those required in Theorem 4.3.

Let $M^m$ be an $m$ dimensional K\"ahler manifold and let $q$ be a
positive integer less than or equal to $m$. A function $\varphi
:M\rightarrow R$ of class $C^2$ is said to be $q$- {\sl
plurisubharmonic} (resp. {\sl strictly } $q$-{\sl
plurisubharmonic}) if, for each point $x\in M$, the trace of the
restriction of the complex Hessian $H(\varphi )$ to any $q$
dimensional complex vector subspace of $T_xM$ is nonnegative
(resp. positive).

Let $D$ be a relatively compact domain in a K\"ahler manifold
$M^m\ $with smooth boundary $\partial D$.  Recall that if $\varrho
$ is a defining function for $D$, the Levi-form $L(\varrho )$ of
$\varrho $ at $x\in \partial D$ is defined as the restriction of
the complex Hessian $H_x(\varrho )$ to the complex subspace
$H_x^{1,0}(\partial D)$ of $T_x^C(\partial D)$. According to
[Si2], we say that $D$ or $\partial D$ is hyper-$q$-convex (resp.
strongly hyper-$q$-convex) if $D$ has a smooth defining function
$\varrho $ such that the eigenvalues $\lambda _1,...,\lambda
_{m-1}$ of the Levi-form $L(\varrho )$ on the holomorphic
distribution $H_x(\partial D)$ at each point $x\in
\partial D$ satisfy $\sum_{i=1}^q\lambda _{j_i}\geq 0$ (resp.
$>0$) for all $1\leq j_i\leq m-1$. The following Lemma is known
(cf. [Si2], [NS]):

\proclaim {Lemma 6.5} Let $D$ be a hyper-$(m-1)$-convex domain in
a K\"ahler manifold $M^m$ with smooth boundary. Let
$f:D\rightarrow N$ be a harmonic map into a K\"ahler manifold with
strongly seminegative curvature such that $f\in C^2(
\overline{D},N)$ and $\overline{\partial }_bf=0$ on $\partial D$.
Then $f$ is pluriharmonic. \endproclaim

We now give an alternative proof of the following result in [CL].

\proclaim{Proposition 6.6} ([CL]) Let $M^m$ be a K\"ahler manifold
and $D$ a hyper-$(m-1)$-convex domain in $M$ with smooth boundary.
Suppose there exists a plurisubharmonic function $\Psi \in
C^2(\overline{D})$ so that $(\Psi _{i \overline{j}})$ has at least
two positive eigenvalues at some point in $D$. Let
$u:\overline{D}\rightarrow N$ be a smooth map into a complete
K\"ahler manifold with strongly seminegative curvature. If $u$
satisfies the tangential Cauchy-Riemann equation
$\overline{\partial }_bu=0$ on $\partial D$, then there exists a
unique holomorphic extension $f$ of $u$.
\endproclaim
\demo{Proof} Clearly the strongly semi-negativity of the curvature
tensor of $N$ implies the nonpositivity of the sectional curvature
of $N$. By the existence theorem of Hamilton [Ha] and Schoen [Sc],
there exists a unique harmonic map $f$ $:D\rightarrow N$ smooth up
to boundary which solves the Dirichlet problem $f=u$ on $\partial
D$.

By Lemma 6.5, we know that $f$ is pluriharmonic. Clearly $\Psi$
satisfies the conditions of Corollary 6.3. Hence we conclude that
$f$ is holomorphic. \qed
\enddemo

\remark{Remark 6.2} The notions of hyperconvex domains in [CL] is
equivalent to the notion of hyper-$(m-1)$-convex domains defined
here.
\endremark

\proclaim{Theorem 6.7} Let $\overline{D}\subset M$ be a compact
connected domain in a K\"ahler manifold with a smooth defining
function $\Psi $. Set $D(c,0)=\{x\in D:c<\Psi (x)<0\}$ for some
$c<0$. Suppose the function $\Psi $ is $(m-1)$-plurisubharmonic on
$D(c,0)$ and strictly $(m-1)$-plurisubharmonic at some point of
$D(c,0)$. Suppose $N$ is a complete K\"ahler manifold with
strongly seminegative curvature. Let $u:\overline{D}\rightarrow N$
be a smooth map satisfying the tangential Cauchy-Riemann equation
$\overline{\partial }_bu=0$ on $\partial D$. Then there exists a
unique holomorphic extension of $u$.
\endproclaim

\demo{Proof} As in the proof of Proposition 6.6, we have a unique
harmonic map $ f:D\rightarrow N$ smooth up to boundary which
solves the Dirichlet problem $ f=u$ on $\partial D$.

Let $\varepsilon _1\leq \varepsilon _2\leq \cdots \leq \varepsilon
_m$ be the eigenvalues of the complex Hessian $H(\Psi )$. Set
$X=\frac 12\nabla \Psi $ and $k_1(x)=\sum_{i=1}^{m-1}\varepsilon
_i(x)$ for $x\in D$. From (1.13) , Lemma 4.9 and Remark 4.3, we
have
$$
\int_{\partial D}S_\sigma (X,\nu )-\int_{\partial D_c}S_\sigma
(X,\nu )\geq \int_{D(c,0)}<S_\sigma ,\nabla \theta _X>
$$
where $D_c=\{x\in D:\Psi (x)\leq c\}$. By carrying out similar
arguments as in the proof of Theorem 6.1, we may deduce that
$$
\int_{D(c,0)}k_1(x)|\overline{\partial }f|^2\leq 0
$$
The assumptions on $\Psi$ mean that $k_1$ is quasipositive on
$D(c,0)$. It follows that $f:D\rightarrow N$ is holomorphic. \qed
\enddemo

Finally we give the following result:

\proclaim{Theorem 6.8} Let $M$ be as in Theorem 4.3 (resp. 4.7).
Suppose $f:M\rightarrow N$ is either a pluriharmonic map into a
K\"ahler manifold or a harmonic map into a K\"ahler manifold with
strongly semi-negative curvature. If
$$
\underset{R\rightarrow \infty }\to{\lim \inf} (R\int_{\partial
B_\Phi (R)}|\overline{
\partial }_bf|^2)=0 \quad
\text{(resp. }\underset{r\rightarrow \infty
}\to{\lim\inf}(r\int_{\partial B_r(x_0)}| \overline{\partial
}_bf|^2)=0\text{)} \tag{6.9}
$$
then $f:M\rightarrow N$ is holomorphic.
\endproclaim
\demo{Proof} Without loss of generality, we assume that $M$
satisfies the conditions in Theorem 4.3. Set $X=\Phi \nabla \Phi$
. For a regular value $ R $ of $\Phi $, we have the unit outward
normal vector field $\nu =\frac{\nabla \Phi }{|\nabla \Phi |}$
along $\partial B_\Phi (R)$. Using (6.4), we derive that
$$
\aligned S_\sigma (X,\nu )&=R|\nabla \Phi |[\frac{|\sigma
|^2}2-<\sigma (\nu ),\sigma
(\nu )>] \\
&\leq R\sqrt{k_2}|\overline{\partial }_bf|^2
\endaligned
\tag{6.10}
$$
Likewise we have (6.2) on each $B_\Phi (R)$. By the assumption,
there exists a sequence $\{R_i\}$ such that
$$
\lim_{i\rightarrow \infty }R_i\int_{\partial
B_{r_i}(x_0)}|\overline{\partial }_bf|^2=0  \tag{6.11}
$$
It follows that
$$
\int_M\lambda |\overline{\partial }f|^2=0
$$
where $\lambda$ is given by (4.4). Therefore $f$ is holomorphic.
\qed
\enddemo

{\bf Acknowledgments}:The author would like to thank Professor
J.G. Cao, Professor N. Mok, Professor M.C. Shaw, Professor S.W.
Wei and Dr. Q.C. Ji for their valuable suggestions and helpful
discussions. He would also like to thank Dr. H.Z. Lin and Dr. G.L.
Yang for their careful reading of the manuscript.

\vskip 0.5 true cm

\vskip 1 true cm

Institute of Mathematics

Fudan University, Shanghai 200433

P.R. China

And

Key Laboratory of Mathematics

for Nonlinear Sciences

Ministry of Education

\vskip 0.2 true cm yxdong\@fudan.edu.cn

\enddocument